\documentclass{article}
\usepackage{graphicx} 
\usepackage{tabularray}
\usepackage{adjustbox}
\usepackage{amssymb}
\usepackage{amsmath}
\UseTblrLibrary{booktabs}
\usepackage{float}
\usepackage{xcolor}
\usepackage{authblk}
\usepackage[margin=1.2in]{geometry}
\usepackage{subcaption}
\usepackage{comment}
\usepackage{amsthm}
\theoremstyle{plain}
\newtheorem{deff}{Definition}

\newtheorem{problem}{Problem}

\newcommand{\R}{\mathbb{R}}
\newcommand{\Z}{\mathbb{Z}}
\newcommand{\C}{\mathbb{C}}

\title{Learning cardiac activation and repolarization times with operator learning}
\author[1]{Edoardo Centofanti\textsuperscript{*}}
\author[2]{Giovanni Ziarelli\textsuperscript{*}}
\author[3]{Nicola Parolini} 
\author[2]{Simone Scacchi}
\author[3]{Marco Verani} 
\author[1]{Luca F. Pavarino} 

\affil[1]{\small Dipartimento di Matematica, Università di Pavia, Via Adolfo Ferrata, 5, Pavia, 27100, Italy}
\affil[2]{\small Dipartimento di Matematica, Università di Milano, Via Cesare Saldini, 50, Milano, 20133, Italy} 
\affil[3]{\small MOX Laboratory, Dipartimento di Matematica, Politecnico di Milano, Via Edoardo Bonardi, 9, Milano, 20133, Italy}

\date{}

\begin{document}

\maketitle
\begingroup
\renewcommand\thefootnote{\textsuperscript{*}}
\footnotetext{The first two authors have contributed equally.}
\endgroup

\begin{abstract}
\noindent Solving partial or ordinary differential equation models in cardiac electrophysiology is a computationally demanding task, particularly when high-resolution meshes are required to capture the complex dynamics of the heart.
Moreover, in clinical applications, it is essential to employ computational tools that provide only relevant information, ensuring clarity and ease of interpretation. 
In this work, we exploit two recently proposed operator learning approaches, namely Fourier Neural Operators (FNO) and Kernel Operator Learning (KOL), to learn the operator mapping the applied stimulus in the physical domain into the activation and repolarization time distributions.
These data-driven methods are evaluated on synthetic 2D and 3D domains, as well as on a physiologically realistic left ventricle geometry.  
Notably, while the learned map between the applied current and activation time has its modelling counterpart in the Eikonal model, no equivalent partial differential equation (PDE) model is known for the map between the applied current and repolarization time.
Our results demonstrate that both FNO and KOL approaches are robust to hyperparameter choices and computationally efficient compared to traditional PDE-based Monodomain models. These findings highlight the potential use of these surrogate operators to accelerate cardiac simulations and facilitate their clinical integration.
\end{abstract}

\section{Introduction}
Computational modeling of cardiac electrophysiology has become a fundamental tool for understanding heart function, diagnosing cardiac conditions, and developing therapeutic interventions. Recent years have witnessed significant advances in both mathematical modeling, numerical techniques and computational capabilities, enabling increasingly sophisticated simulations of cardiac electrical activity. Despite these advances, the computational complexity of high-fidelity cardiac models remains a substantial challenge, particularly for large-scale simulations, real-time applications, and clinical decision support systems. The Bidomain model~\cite{cla11, col14, kee09, tra24, sun07} serves as the gold standard for describing the propagation of extra- and intracellular potentials in cardiac tissue, though its computational complexity presents significant challenges for large-scale simulations. In several works, in particular involving electromechanical coupling or interactions with fluid dynamics models of the human heart, researchers often turn to the more computationally efficient Monodomain model~\cite{col14, kee09} as an alternative. The latter emerges as a simplified version of the Bidomain model where the intra- and extracellular conductivities are proportional. This simplification results in a model that maintains reasonable accuracy while significantly reducing the computational demand. The Monodomain model is constituted by a system of partial differential equations that describe the spatiotemporal evolution of the transmembrane potential and associated gating or recovery variables, which either represent the probability of ionic species flowing through the membrane or serve as recovery variables designed to reproduce observed phenomenological potentials (see, \textit{e.g.}, the two-variable models derived from the FitzHugh-Nagumo model~\cite{fit61}). A more computationally efficient approach is provided by Eikonal models~\cite{col14, kee09}, which focus on the evolution of cellular excitation wavefronts rather than the complete spatial and temporal reconstruction of ionic action potentials: these models are extremely computationally cheap, though they provide less details regarding the upstroke thin layer of propagation.

\begin{figure}[t]
    \centering
    \includegraphics[width=\textwidth]{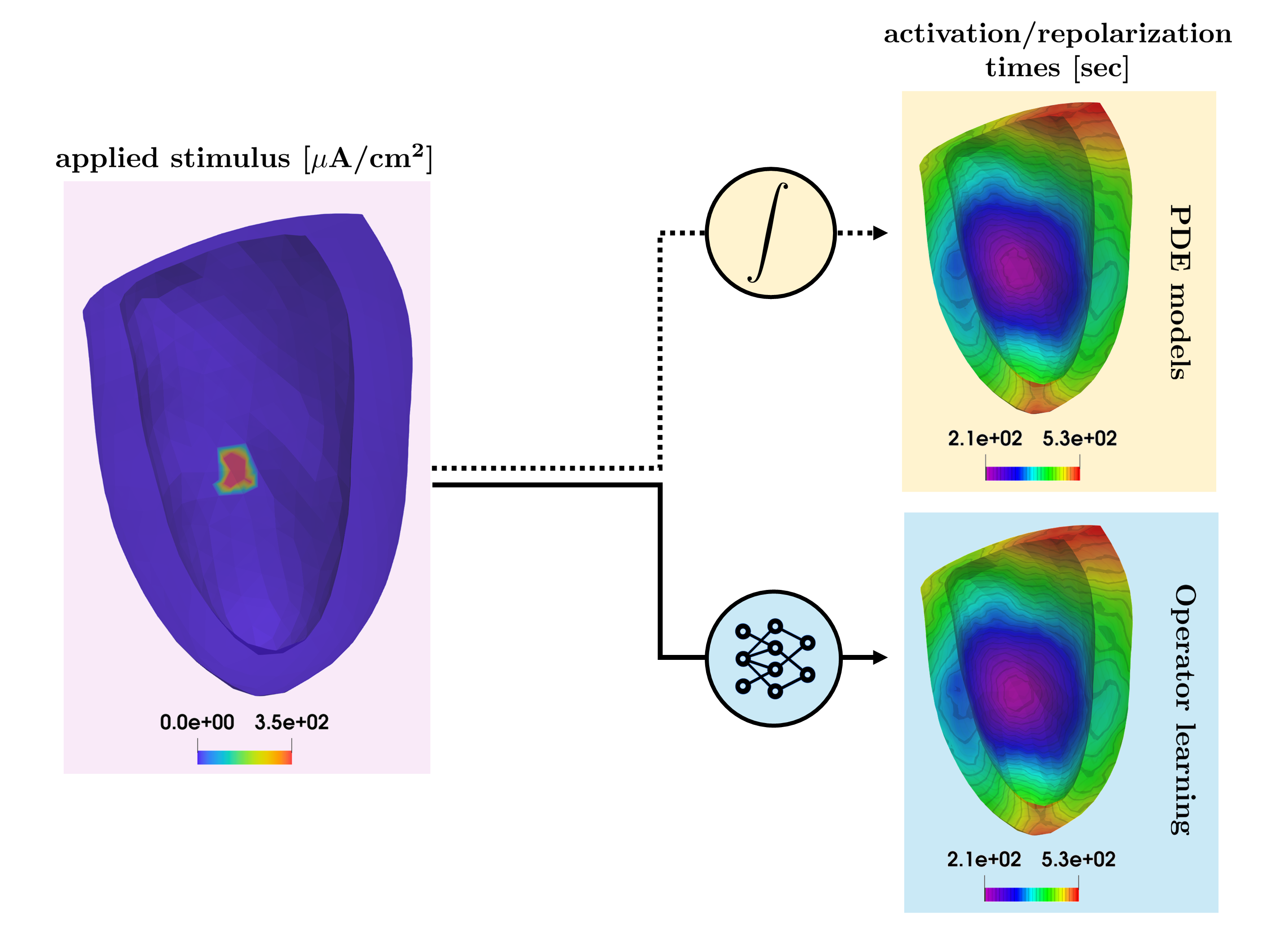}
    \caption{Schematic representation of the EP problem to address. In particular, we aim at reconstructing the activation/repolarization times of the cardiac tissue given the initial stimulus applied for 1 ms.}
    \label{fig:introscheme}
\end{figure}

Despite the extensive range of mathematical models and numerical schemes available for addressing electrophysiological (EP) problems, their computational burden remains a significant concern.
Furthermore, the primary interest in solving these models often lies in extracting key informative quantities that can significantly assist clinicians, such as activation and repolarization times within the cardiac domain.
The activation and repolarization times serve as critical markers in cardiac electrophysiology, providing essential information about the heart's electrical function: the activation time refers to the time when cardiac cells begin their depolarization process, whilst repolarization time denotes the time when cells return to their resting state.
These markers are fundamental for understanding cardiac conduction patterns, identifying arrhythmogenic substrates~\cite{cor09}, and evaluating the effects of drugs or interventions~\cite{var11, and02}. Traditional approaches for computing these times typically require solving the full Monodomain or Bidomain models and extracting the times from the resulting action potential waveforms, a process that is computationally intensive. Moreover, while activation times can be evaluated more efficiently through the Eikonal model, repolarization times lack of a classical Eikonal-like counterpart.

The emergence of scientific machine learning offers new opportunities to address these computational challenges~\cite{cen24b, fre20, kis20, reg20, tra21, ten22, reg24}. In particular one of its main branches, operator learning, aims to approximate unknown operators that map between potentially infinite-dimensional functional spaces. Given pairs of input/output functional data $(u, f)$, where $u \in \mathcal{U}$ and $f \in \mathcal{V}$ are functions defined on domains $\Omega$ and $\Omega'$ respectively, the goal is to learn an approximation of an operator $\mathcal{G}: \mathcal{U}\to \mathcal{V}$ using machine learning architectures. Among the various recently-proposed operator learning architectures (see, \textit{e.g}, \cite{lu19, go22, tri23}), Fourier Neural Operators (FNOs)~\cite{li20} have emerged as a powerful approach based on the Neural Operator paradigm~\cite{kov21b}, parameterizing the integral kernel layers within the architecture in the Fourier space and allowing efficient learning of mappings between function spaces with resolution independence. FNOs have shown comparable performances for equispaced domains with respect to the vanilla Deep Operator Networks \cite{lu22}. Another promising approach is Kernel Operator Learning (KOL)~\cite{bat24}, which builds on standard kernel regression arguments to approximate the mapping between function spaces. Compared to other neural operator methodologies, the key advantage of the Kernel Operator Learning (KOL) approach lies in its non-iterative formulation; the operator is obtained by solving a (potentially large) symmetric and positive definite linear system, thereby eliminating the need for iterative training procedures typically required in neural network-based frameworks.

Using operator learning techniques to predict activation and repolarization times based on inputs such as tissue conductivity, fiber orientation, and stimulus location can help reduce computational bottlenecks in traditional cardiac modeling.
These approaches will also improve efficiency, making computational tools more accessible for research and clinical practice.
In this work, we take a step in this direction by learning the mapping between an initially applied current stimulus and the activation/repolarization times at each physical point in the considered 2D or 3D domains, as schematically represented in Figure~\ref{fig:introscheme}.
Specifically, we compare the performances of FNO e KOL in terms of training time, testing time, memory usage and accuracy in testing, and we assess the potentialities of both strategies for retrieving fast and accurate simulations.

This paper is structured as follows. In Section~\ref{sec:methods} we review the mathematical models (Section~\ref{sec:pdemodels}) and we introduce and formalize FNO and KOL for the problem at stake and used in the numerical experiments (Sections~\ref{sec:ol}-\ref{sec:kol}). In Section~\ref{sec:numerical} we detail the dataset generation (Section~\ref{sec:datasetgen}) and presents numerical results for both 2D and 3D cases, where we discuss the computational performance of the trained models (Sections~\ref{sec:2Dcase}-\ref{sec:3dunstcase}).
Finally, in Section~\ref{sec:conclusions} we draw some concluding remarks. 

\section{Methods}\label{sec:methods}
\subsection{PDE models}\label{sec:pdemodels}
Mathematical models of electrophysiology (EP) play a crucial role in understanding and simulating the electrical activity of cardiac tissue. These models describe the evolution of the transmembrane potential and ionic currents, catching the fundamental mechanisms of excitation and propagation. There is a plethora of EP models that may be used depending on the required level of detail (see, \textit{e.g.}, \cite{col14}), ranging from ionic models, based on systems of ODEs, to other macroscopic formulations.

The Bidomain model \cite{cen24} describes the propagation of the extra- and intracellular potentials and it is widely used in order to model the electrical activation in the myocardium. However, as already mentioned, its computational complexity makes it challenging for large-scale simulations \cite{afr22, buc24} and the Monodomain model \cite{ver16} is often preferred as a more efficient alternative. 

The Monodomain model is derived from the cardiac Bidomain model when the intra- and the extracellular conductivities, $\mathbf{D}_i$ and $\mathbf{D}_e$ respectively, satisfy the following relationship \cite{col14}:

\begin{equation}
    \mathbf{D}_e = \lambda \mathbf{D}_i
    \label{eq:assumpCond}
\end{equation}

where $\lambda \in \R$ is a constant. The model reads as follows:

\begin{equation}
\left\{
\begin{array}{ll}
\displaystyle \chi C_m \frac{\partial v}{\partial t} - \frac{\lambda}{1 + \lambda} \text{div}(\mathbf{D}_i \nabla v) + I_{\text{ion}}(v, \mathbf{w}) = I_{\text{app}} & \text{in } \Omega \times (0, T),\vspace{0.2cm}\\
\displaystyle \frac{\partial \mathbf{w}}{\partial t} - \mathbf{R}(v, \mathbf{w}) = 0 & \text{in } \Omega \times (0, T),\vspace{0.2cm} \\
\displaystyle \frac{\partial \mathbf{\mathbf{c}}}{\partial t} - \mathbf{C}(v, \mathbf{w},\mathbf{c}) = 0 & \text{in } \Omega \times (0, T),\vspace{0.2cm}\\
\displaystyle \mathbf{n}^\top \mathbf{D}_i \nabla v = 0 & \text{on } \partial \Omega \times (0, T),\vspace{0.2cm}\\
\displaystyle v(\mathbf{x}, 0) = v_0(\mathbf{x}), \quad \mathbf{w}(\mathbf{x}, 0) = \mathbf{w}_0(\mathbf{x}) & \text{in } \Omega.
\end{array}
\right.
\label{eq:monodomain}
\end{equation}

where $v$ is the transmembrane potential, which represents the difference between the intra- and extracellular potentials, $C_m$ and $\chi$ are the membrane capacitance per unit area and the membrane surface area per unit volume respectively, $\lambda$ is the constant in~\eqref{eq:assumpCond}, $I_{\text{ion}}$ is a ionic current density representing the flow of ionic species through the cellular membrane and $I_{\text{app}}$ is an applied stimulus in time. This latter term depends on $v$, as well as on the gating or recovery variables of the ionic model, $\mathbf{w}$. These variables either describe the probability of ionic species flowing through the membrane or serve as recovery variables designed to reproduce an observed phenomenological potential, as seen in two-variable models derived from the FitzHugh-Nagumo \cite{izh06} model. They are coupled to the reaction diffusion PDE through a system of differential equations describing their evolution in time as well as the dynamics of the ionic concentrations $\mathbf{c}$, which are ruled by often nonlinear functions $\mathbf{R}(v,\mathbf{w})$ and $\mathbf{C}(v,\mathbf{w},\mathbf{c})$.
If there is no injection of current in the extracellular space, \eqref{eq:monodomain} can be considered as a good approximation of the Bidomain model.

An alternative approach for modelling the evolution of the cellular excitation wavefront is based on Eikonal models.
Starting from the Bidomain model combined with a simplified 
representation of the ionic currents, where we exclude any gating variables and we focus solely on the depolarization phase, the Eikonal-diffusion equation is (see, \textit{e.g.}, \cite{col90}):
\begin{equation}
\left\{
\begin{array}{ll}
\displaystyle c_0 \sqrt{\nabla \psi \, \cdot \, \mathbf{M} \nabla \psi} \, - \, \nabla \cdot (\mathbf{M}\nabla \psi) = 1, & \mathrm{in} \, \Omega,\\[5pt]
\mathbf{M} \nabla \psi \cdot \mathbf{n} = 0, & \mathrm{on} \, \partial \Omega_n, \\[5pt]
\psi = \psi_D, & \mathrm{on} \, \partial \Omega_D,
\end{array}
\right.
\label{eq:eikonal}
\end{equation}

where $\psi$ is the unknown activation time of cardiac cells, \textit{i.e.} the instant when the transmembrane potential first crosses a predefined threshold during an action potential, marking the onset of electrical excitation.
Moreover, $c_0$ represents the estimated velocity of the depolarization wave along the fiber direction of the cardiac tissue for planar wavefronts, and $\mathbf{M} = \frac{1}{\chi C_m}\frac{\lambda}{1 + \lambda} \mathbf{D}_i$. Well-posedness is ensured by imposing the Dirichlet boundary condition on a non-trivial portion of the domain boundary.
An alternative model, known as Eikonal-curvature, has been proposed in \cite{ken91}, and reads as
\begin{equation}
    c_0 \sqrt{\nabla \psi \, \cdot \, \mathbf{M} \nabla \psi} \, - \, \sqrt{\nabla \psi \, \cdot \, \mathbf{M} \nabla \psi} \, \nabla \cdot \bigg (\dfrac{\mathbf{M}\nabla \psi}{ \sqrt{\nabla \psi \, \cdot \, \mathbf{M} \nabla \psi}} \bigg ) = 1,
\end{equation}
endowed with the same boundary conditions of \eqref{eq:eikonal}.
In this model the diffusion term is proportional to the anisotropic generalized mean curvature~\cite{col14, kee09}.
In both equations, contour levels of $\psi$ display the activation times at each prescribed location.
Eikonal models are particularly attractive from a computational perspective, since they involve a single steady-state PDE that, despite being nonlinear, does not require coupling with ODE systems. More importantly, unlike the transmembrane potential, the activation time lacks internal or boundary layers, eliminating the need for special mesh restrictions.
Finally, we remark that activation and repolarization times can be both derived by post-processing the solution of the Monodomain model. However, while activation times can be obtained by solving Eikonal equations, no Eikonal-based formulation has been proposed in the literature to directly extract repolarization times.

In the next section we are going to introduce the operator learning tools that will enable us to construct surrogate models for both activation and repolarization times.

\subsection{Operator Learning: Basic Principles}\label{sec:ol}

Operator Learning (OL) aims to approximate an unknown operator 

\begin{equation*}
    \mathcal{G}: \mathcal{A} \to \mathcal{U},
\end{equation*}

which maps between two infinite-dimensional functional spaces $\mathcal{A}$ and $\mathcal{U}$. Given data pairs $(a, u)$, where $a \in \mathcal{A}$ and $u \in \mathcal{U}$ are functions defined on bounded domains $\Omega\subset\mathbb{R}^n$, $\Omega'\subset\mathbb{R}^m$ and sampled over $\mathbb{R}^n$ and $\mathbb{R}^m$ respectively, the goal is to learn an approximation of $\mathcal{G}$ by mean of machine learning surrogate model.
The problem can be formalized as follows:

\begin{problem}\label{prob:setting}
    Let us consider $\{a_i,u_i\}_{i=1}^N$ samples in $\mathcal{A}\times\mathcal{U}$, such that
    \begin{equation}
        \mathcal{G}(a_i) = u_i, \qquad \text{with } i=1,2,\dots, N.
    \end{equation}
    We define the observation operators $\phi : \mathcal{A}\to \mathbb{R}^n$ and $\varphi : \mathcal{U}\to \mathbb{R}^m$ acting on the input and the output functions, respectively. The aim of operator learning is approximating the operator $\mathcal{G}$ through the observation of input/output pairs $\{\phi(a_i),\varphi(u_i)\}$.
\end{problem}
In this work we consider a natural choice for the observation operators $\phi$ and $\varphi$, namely the pointwise evaluation at specific collocation points $\{\mathbf{x}_k\}_{k=1}^n \in \Omega$ and $\{\tilde{\mathbf{x}}_k\}_{k=1}^m \in \Omega'$ respectively, which in general can be different, \textit{i.e.}
\begin{equation}
  \begin{split}
  &\phi : a \rightarrow A := (a({\mathbf{x}}_1), a({\mathbf{x}}_2) \hdots a({\mathbf{x}}_n))^T \in \mathbb{R}^{n \cdot d_a},\\
  &\varphi : u \rightarrow U := (u(\tilde{\mathbf{x}}_1), u(\tilde{\mathbf{x}}_2) \hdots u(\tilde{\mathbf{x}}_m))^T \in \mathbb{R}^{m \cdot d_u}.
  \end{split}
\end{equation}
However, for simplicity, in this work we consider the same collocation points, namely $n=m$ and $\{\mathbf{x}_k\}_{k=1}^n = \{\tilde{\mathbf{x}}_k\}_{k=1}^n$.
We work in a supervised learning framework assuming we are given the training dataset $\{ A_i, U_i\}_{i=1}^N$ where, consistently with our notation, 
$$A_i := (a_i({\mathbf{x}}_1), a_i({\mathbf{x}}_2) \hdots a_i({\mathbf{x}}_n))^T \in \mathbb{R}^{n \cdot d_a} \; \mathrm{and} \; U_i := (u_i(\mathbf{x}_1), u_i(\mathbf{x}_2) \hdots u_i(\mathbf{x}_n))^T \in \mathbb{R}^{n \cdot d_u}.$$
Among all the possible operators ranging between $\mathcal{A}$ and $\mathcal{U}$ we aim at finding the one which minimizes the error in a prescribed suitable norm and assuming tailored paradigms for the approximated operator, \textit{e.g.} machine learning architectures.
In the following, we briefly expand on the two methods employed for the reconstruction of activation and repolarization times, namely Fourier Neural Operators and Kernel Operator Learning.

\subsection{Fourier Neural Operators}\label{sec:fno}

Fourier Neural Operators (FNOs) \cite{li20} are operator learning schemes based on the Neural Operator paradigm \cite{kov21b}. 

\begin{figure}[t]
    \centering
    \includegraphics[width=0.8\textwidth]{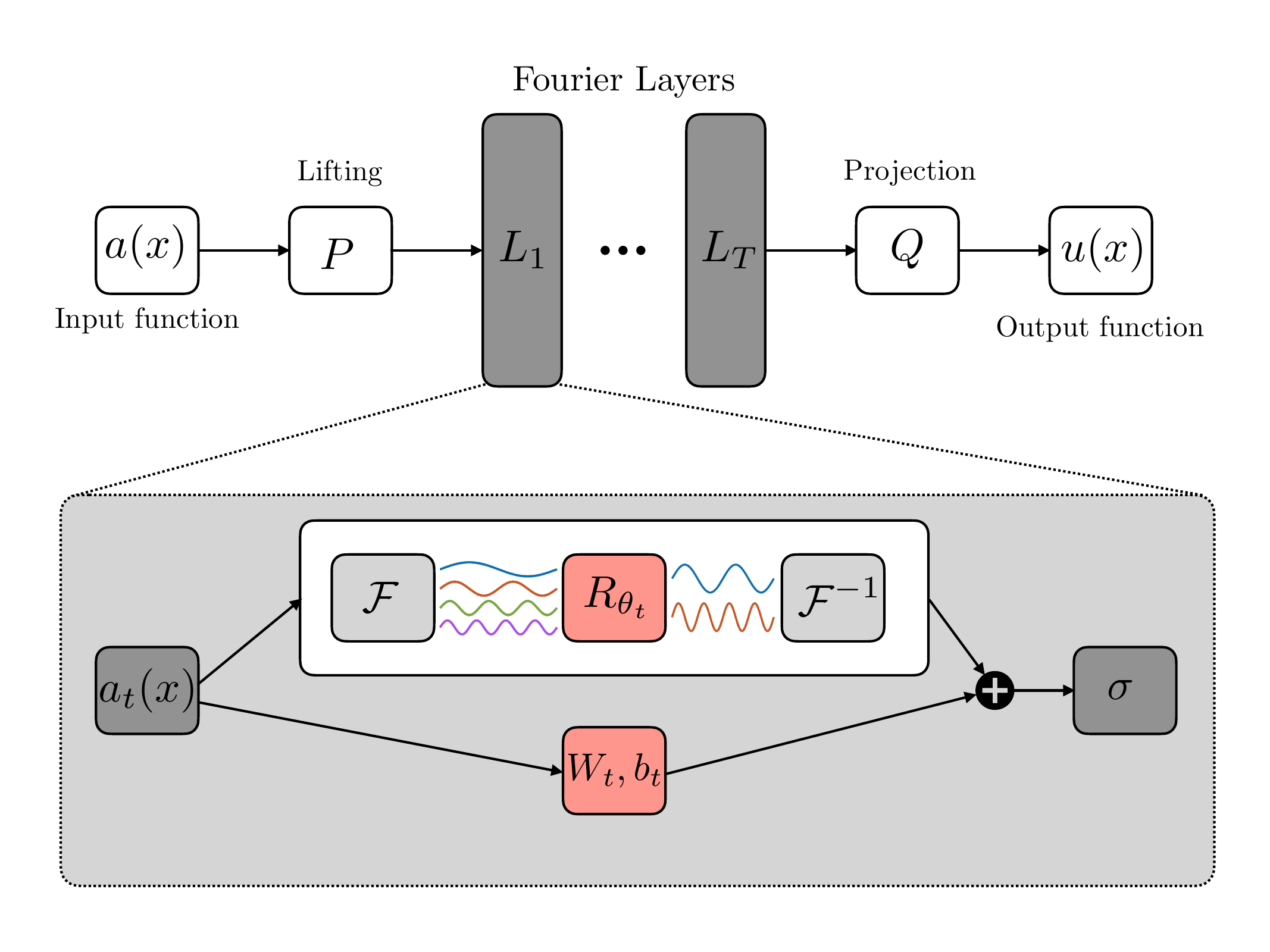}
    \caption{Schematic architecture of Fourier Neural Operators (FNO).}
    \label{fig:fno_arch}
\end{figure}

\begin{deff}[Neural Operator]
    Assuming the setting introduced in Problem~\ref{prob:setting}, we define the Neural Operator as the architecture $\hat{\mathcal{G}}_\theta:\R^n\rightarrow\R^m$ with the following structure:
    
    \begin{equation}\label{eq:neural_op}
    \hat{\mathcal{G}}_\theta := Q\circ \sigma_T\left( W_{T-1}+K_{T-1} + b_{T-1} \right)\circ \cdots\circ \sigma_1\left( W_0+K_0 + b_0 \right)\circ P
    \end{equation}

where each inner operator represents a layer which maps the $t$-th hidden representation $a_t$ to the next one $a_{t+1}$, for $t=0\dots T-1$, following the scheme reported in Figure~\ref{fig:fno_arch}. More in detail, a Neural Operator is composed as follows:

\begin{itemize}
    \item $P:\R^n\rightarrow\R^{d_{a_{0}}}$, a {lifting} operator, namely a pointwise function mapping the observed input to its first hidden representation, \textit{i.e.} $\{\phi(u): \mathcal{U}\rightarrow \R^{n}\}\mapsto\{a_0: \Omega_0\subset\R^{d_0}\rightarrow\R^{d_{a_0}}\}$, with $d_{a_0} > d_0$. This operation is performed by a fully local operator, and its action is pointwise, namely $(P(u))(x) = P(u(x))$ for $x\in \Omega$;

    \item A composition of operators $\{a_t:\Omega_t\rightarrow \R^{d_{a_t}}\} \rightarrow \{a_{t+1}:\Omega_{t+1}\rightarrow \R^{d_{a_{t+1}}}\}$, for $t=0,\dots,T-1$, each one defined as the sum of a local linear operator $W_t \in \R^{d_{a_{t+1}} \times d_{a_t}}$, a non-local integral kernel operator $K_t$, and a bias function $b_t:\Omega_{t+1} \rightarrow \R^{d_{a_{t+1}}}$. This sum is then composed with a fixed pointwise nonlinearity $\sigma_t$, called activation function, namely $(\sigma_t(a_{t+1}))(x) = \sigma(a_{t+1}(x))$ for any $x \in \Omega_{t+1}$. We set $\Omega_0 = \Omega$ and $\Omega_T = \Omega'$, and assume that each $\Omega_t \subset \R^{d_t}$ is a bounded domain.

    
    \item $Q:\R^{d_{a_{T}}}\rightarrow\R^{m}$, a {projection} operator, namely a pointwise function, mapping the last hidden representation $\{a_T:\Omega'\rightarrow\R^{d_{a_{T}}}\}\mapsto\{\varphi(u):\mathcal{U}\rightarrow\R^{m}\}$ to the observed output function. Since $d_{a_{t}}>m$, this is a projection step performed by a fully local operator, namely $(Q(a_T))(x) = Q(a_T(x))$ for every $x\in \Omega'$.
\end{itemize}
\end{deff}

Input and output dimensions $d_{a_{0}},\cdots,d_{a_{T}}$ as well as the domains of definition $\Omega_1,\cdots, \Omega_{T-1}$ are hyperparameters of the architecture.

In this work, we follow the first definition for the integral kernel operators proposed in \cite{kov21b}. As a technical remark, we denote as $C(\Omega_{t+1}\times \Omega_t,\, \R^{d_{a_{t+1}}\times d_{a_{t}}})$ the Banach space of continuous functions from $ \Omega_{t+1} \times \Omega_t$ to $\R^{d_{a_{t+1}}\times d_{a_{t}}}$ equipped with the sup-norm $||g||_\infty := \sup \{ |g(\tau)|:\tau\in\Omega_{t+1}\times\Omega_t \},$ $g\in C(\Omega_{t+1}\times \Omega_t,\, \R^{d_{a_{t+1}}\times d_{a_{t}}})$.

In this context, given $\kappa^{(t)}\in$ $C(\Omega_{t+1}\times \Omega_t; \R^{d_{a_{t+1}}\times d_{a_{t}}})$ such that $\kappa^{(t)}(x,y) = \kappa^{(t)}(x-y)$, for all $(x,y)\in\Omega_{t+1}\times \Omega_t$, the integral kernel operator $K_t$ is defined as
\begin{equation}\label{eq:kappat}
    (K_t(a_{t}))(x) = \int_{\Omega_t} \kappa^{(t)}(x-y) a_{t}(y) \ dy \quad \forall x \in \Omega.
\end{equation}

At this point, we can thus make precise the single hidden layer update rule as
\begin{equation}
    a_{t+1}(x) = \sigma_{t+1} \left( W_ta_{t}\left( x \right) + \int_{\Omega_t} \kappa^{(t)}(x,y)a_{t}(y)\,d\nu_t(y) + b_t(x)\right)\qquad \forall x\in \Omega_{t+1}.
\end{equation}
From this general structure for the Neural Operator, one can derive the FNO by discretizing the integral kernel (\ref{eq:kappat}) in the Fourier space.
The key feature of the FNO is that the integral kernel is parameterized in the Fourier space. Therefore, in the following we recall the Fourier transform $\mathcal{F}:L^2(\Omega,\mathbb{C}^n)\to \ell^2(\mathbb{Z}^d,\mathbb{C}^n)$ and anti-transform $\mathcal{F}^{-1}$. Given $a\in L^2(\Omega,\mathbb{C}^n)$ and $\hat{a}\in\ell^2(\mathbb{Z}^d,\mathbb{C}^n)$, we define
\begin{equation}
    \begin{split}
        & \left( \mathcal{F}a \right)_j(k) =\hat{a}_j(k) = \langle a_j, \varphi_k\rangle_{L^2(\Omega,\mathbb{C}^n)},\qquad j\in\{1,\cdots,n\},\quad k\in\mathbb{Z}^d, \\
        &\left( \mathcal{F}^{-1} \hat{a} \right)_j(x) = \sum_{k\in\mathbb{Z}^d}\hat{a}_j(k)\varphi^{-1}_k(x),\qquad j\in\{1,\cdots,n\},\quad x\in\Omega \\
        &\varphi_k(x) := e^{-2\pi ik \cdot x},\qquad x\in \Omega.
    \end{split}
\end{equation}
For the FNO, the domain considered for each layer is the periodic torus $\Omega_t = \mathbb{T}^d = [0,2\pi]^d$, although for a general input it is sufficient to take its periodic extension rescaled in $[0,2\pi]^d$. 

The integral kernel $\kappa^{(t)}$ in~\eqref{eq:kappat} is a function parameterized by some parameters $\theta_t$ belonging to a suitable space $\Theta_{t}\subset\R^{d_{a_{t+1}}\times d_{a_{t}}}$. Thus we write $\kappa^{(t)}=\kappa^{(t)}_{\theta_t}$. From~\eqref{eq:kappat}, we set $\kappa^{(t)}_{\theta_t}(x, y) = \kappa^{(t)}_{\theta_t}(x-y)$ and we apply the convolution theorem for the Fourier transform,
\begin{equation}\label{eq:FNOKernel}
    (K_t(a_{t}))(x) = \int_{\Omega} \kappa^{(t)}_{\theta_t}(x-y)a_{t}(y)\,dy = \mathcal{F}^{-1}\left( \mathcal{F}\left(\kappa^{(t)}_{\theta_t}\right)\cdot\mathcal{F}(a_{t}) \right)(x),\qquad \forall x\in\Omega_t,
\end{equation}
\noindent
where the dot operation is defined as
\begin{equation} \label{eq:hadamard}
    \mathcal{F}\left(\kappa^{(t)}_{\theta_t}\right) \cdot\mathcal{F}(a_{t}) := \sum_{k \in \Z} \mathcal{F}\left(\kappa^{(t)}_{\theta_t}\right)(k)\cdot \mathcal{F}(a_{t})(k) = \sum_{k \in \Z} \hat{\kappa}^{(t)}_{\theta_t}(k) \cdot \hat{v}_t(k).
\end{equation}
\noindent Considering $k\in\mathbb{Z}^d$ fixed as frequency mode, we have that $\mathcal{F}\left(\kappa^{(t)}_{\theta_t}\right)(k)\in\mathbb{C}^{d_{a_{t+1}}\times d_{a_{t}}}$ and $\mathcal{F}(a_{t})(k)\in\mathbb{C}^{d_{a_{t}}}$. Hence, $\kappa^{(t)}_{\theta_t}$ can be parameterized directly by its Fourier coefficients and~\eqref{eq:FNOKernel} reads as:
\begin{equation}\label{eq:KernelFNO}
   (K_{t}(a_{t}))(x) = \mathcal{F}^{-1}\left( R_{\theta_t}\cdot\mathcal{F}(a_{t}) \right)(x),\qquad\forall x \in \Omega_t,
\end{equation} 
where $R_{\theta_t} = \hat{\kappa}^{(t)}_{\theta_t}(k)$ for $k$ fixed. Fig. \ref{fig:fno_arch} represents the schematic FNO architecture considered in this work.

Numerically, the FNO needs to be approximated from continuous space to discrete ones for dealing with finite-dimensional parameterizations. We call pseudo$(\Psi)$-Fourier Neural Operator ($\Psi$-FNO) the approximated architecture, following \cite{kov21a}. In this case, the domain $\Omega_t$ is discretized with $J\in\mathbb{N}$ points, and therefore $a_{t}$ can be treated as a tensor in $\mathbb{C}^{J\times d_{a_{t}}}$. Since integrals cannot be calculated exactly, this leads to considering the Fourier series truncated at a maximal mode $k_{\max}$, such that
\begin{equation*}    
    k_{\max} = \big|\{k\in\mathbb{Z}^d \ :\ |k_j|\le k_{\max,j},\text{ for } j=1,\dots,d\}\big|.
\end{equation*}
\noindent
In practical implementations, the Fourier transform is replaced by the Fast Fourier Transform (FFT), and the weight tensor $R_{\theta_t}$ is parameterized as a complex-valued tensor $ R_{\theta_t} \in \C^{k_{\max}\times d_{a_{t+1}}\times d_{a_{t}}}$ for $i = 1, \dots, L$. For simplicity, we denote $R_{\theta_t}$ with $R$. 
Finally, we note that if $a_{t}$ is real-valued, we can enforce conjugate symmetry in tensor $R$ for imposing $a_{t+1}$ to be real-valued, namely $R(-k)_{j,l}=\overline{R}(k)_{j,l}$.

\subsection{Kernel Operator Learning}\label{sec:kol}
\begin{figure}[t]
    \centering
    \includegraphics[width=0.4\textwidth]{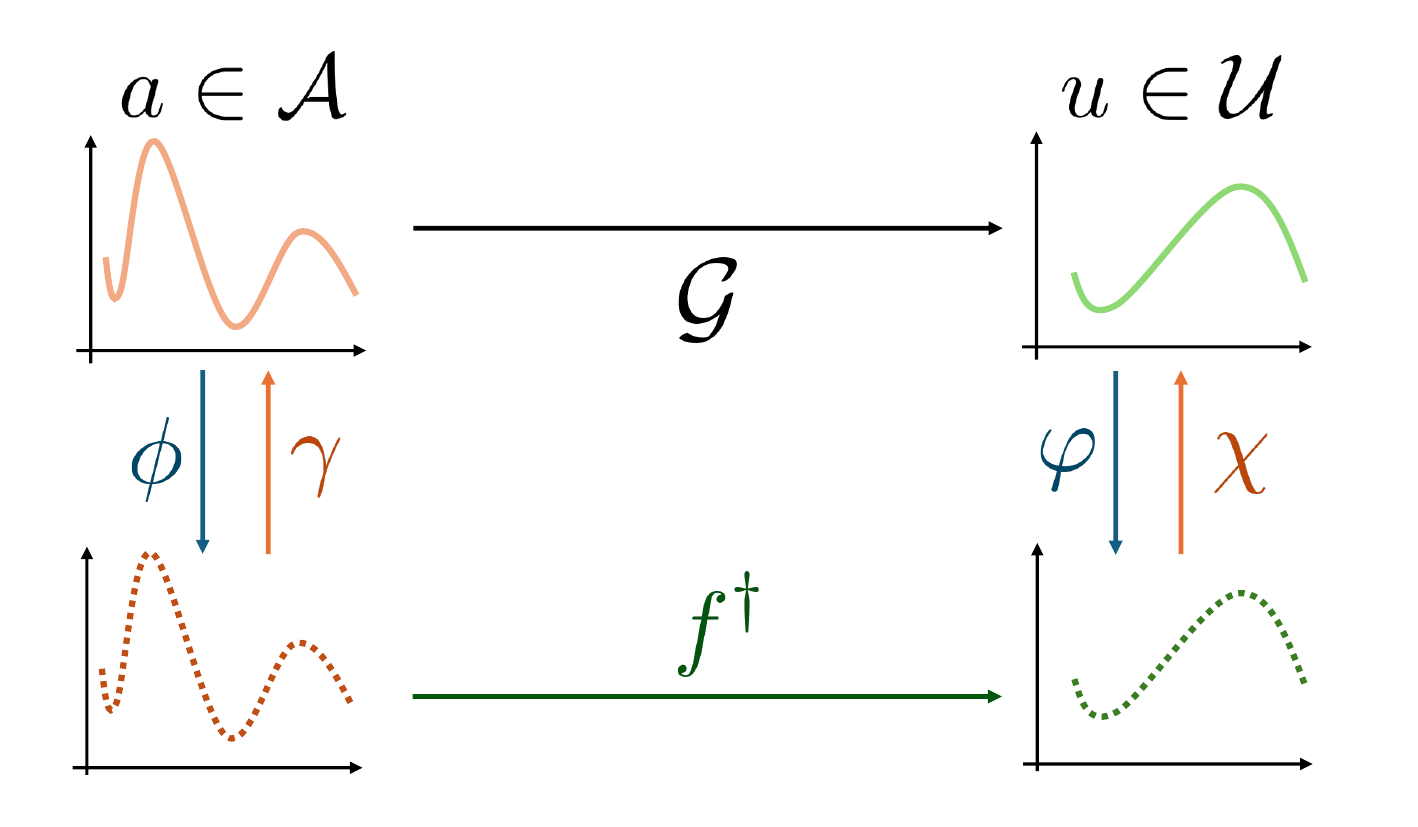}
    \caption{Kernel Operator Learning (KOL) diagram. Starting from the input function $a$, $A$ collects observations of the function at different collocation points through $\phi$. Then, the vector-valued $f^{\dagger}$ processes observations of the input into observations of the output $A = \phi(a)$. Finally, the reconstruction operator is applied to determine the output function $u = \chi(U)$.}
    \label{fig:ol_arch}
\end{figure}
Kernel Operator Learning (KOL) is a recently formalized operator learning technique \cite{bat24}, based on standard kernel regression arguments.
Following the diagram in Figure \ref{fig:ol_arch}, retrieving the approximated operator $\bar{\mathcal{G}}$ is equivalent to determining a vector-valued function $f^{\dagger}: \mathbb{R}^{n \cdot d_a} \rightarrow \mathbb{R}^{n \cdot d_u}$ such that $\bar{\mathcal{G}} = \chi \circ f^{\dagger} \circ \phi$, where $\chi:\R^{n\cdot d_u}\rightarrow \mathcal{U}$ is a reconstruction operator and $\phi$ has been defined in Problem~\ref{prob:setting}.
In the following, we formalize KOL in the scalar-function case, \textit{i.e.} assuming $d_a = d_u = 1$. More specifically, analogously to \cite{zia24}, endowing $\mathcal{A}$ and $\mathcal{U}$ with a Reproducing Kernel Hilbert Space (RKHS) structure and using kernel regression to identify the maps $\gamma$ and $\chi$, the approximated operator can be written explicitly in closed form as
\begin{equation}
  \bar{\mathcal{G}}(a)(\mathbf{x}) = K(\mathbf{x}, \mathcal{X}) K(\mathcal{X},\mathcal{X})^{-1} \left ( \sum_{j=1}^N S(\phi(a), A_j) \alpha_j \right ),
  \label{eq:KOLGeneral}
\end{equation}
where $K$ is the kernel function induced by the RKHS structure of $\mathcal{U}$, the vector $$\mathcal{X} = [\mathbf{x}_1, \mathbf{x}_2 \hdots \mathbf{x}_n]^T \in \mathbb{R}^{n \cdot d}$$ contains the collocation points, and $S:\mathbb{R}^n \times \mathbb{R}^n \rightarrow \mathbb{R}$ is a properly chosen vector-valued kernel.
Moreover, $K(\cdot, \mathcal{X}) : \Omega \rightarrow \mathbb{R}^n$ is a row vector such that $K(\mathbf{x},\mathcal{X})_i = K(\mathbf{x}, \mathbf{x}_i)$, and $K(\mathcal{X},\mathcal{X})$ is an $n \times n$ matrix such that $K(\mathcal{X},\mathcal{X})_{ij} = K(\mathbf{x}_i, \mathbf{x}_j)$.
Since we deal with pointwise observation functions $\phi$ and $\varphi$, $K(\textbf{x}, \mathcal{X})$ computes the evaluation of the linear interpolant of points in $\mathcal{X}$ at \textbf{x}, whilst $K(\mathbf{x}_i, \mathbf{x}_j) = \delta_i^j$. 
Parameters $\{\alpha_j\}_{j=1}^N \in \mathbb{R}^n$ are the kernel regression parameters over the input/output training pairs.

From a computational standpoint, after selecting the scalar kernel $S(\cdot, \cdot)$  for the discrete vector spaces, the problem reduces to solving $n$ linear systems of size $N$ in order to determine the different components of each $\{\alpha_j\}_j$: with a slight abuse of notation we call this matrix $S(A,A) \in \mathbb{R}^{N \times N}$ where each entry is $[S(A,A)]_{ij} = S(A_i, A_j)$.
To achieve this, we employ the Cholesky factorization of the matrix and solve the resulting systems using standard substitution methods.
Moreover, a regularization term is introduced in the regression formulation, with a penalty parameter set to $10^{-10}$.

A key factor influencing the approximation and generalization properties of KOL methods is the selection of the scalar kernel $S$. The optimal kernel function for kernel regression remains a subject of ongoing debate and is application-dependent. For example, some novel approaches involve learning kernels by simulating data-driven dynamical systems, enhancing the scalability of Kernel Regression \cite{owh19}. In this context, we present a few choices for $S$ used in this study, with additional options available in \cite{zia24}.
We note that the proposed kernel functions are specifically tailored for input functions that represent stimuli, such as indicator functions. When working with smoother functions, it is necessary to appropriately adjust the kernel functions based on the specific case at hand.

\begin{itemize}
  \item \textbf{Radial Basis Functions (RBF) kernel}: $S(A_1, A_2) = e^{-\frac{\| \mathbf{c}_{A_1} -  \, \mathbf{c}_{A_2} \|_2^2}{2 \sigma^2}}$, where $ \mathbf{c}_{A_*}$ represents the vector of coordinates of the centroid of non-zero elements of discrete observations in $A_*$.
  It has the interpretation of a similarity measure and it decreases as long as the distance between points increases. In the case of RBF kernel, Kernel Operator Learning has an explicit connection with Gaussian Processes (GP) for regression tasks (see, \textit{e.g.}, \cite{smo98}).
  \item \textbf{Neural Tangent Kernel (NTK)}: Given a neural network regressor $f(x; \theta)$ of depth $d_{\mathrm{nn}}$, width $l_{\mathrm{nn}}$ and activation function $\sigma_{\mathrm{nn}}$, with $\theta_{\mathrm{nn}}$ denoting the vector collecting all weights and biases, we define the family of finite-width Neural Tangent Kernels $ \{ S \}_{\tau > 0}: \mathbb{R}^n \times \mathbb{R}^n \rightarrow \mathbb{R} $ as
  \begin{equation}
    S_\tau(A_i,A_j) := \langle \partial_{\theta_{\mathrm{nn}}} f(A_j; \theta_{\mathrm{nn}}(\tau)) ,\, \partial_{\theta_{\mathrm{nn}}} f(A_i; \theta_{\mathrm{nn}}(\tau)) \rangle,
  \end{equation}
  where $\tau$ represents a fictitious iteration time. 
  It has been proven that, if the initialization of the weights follows the so-called NTK initialization \cite{jac18}, in the infinite-width limit each element in the sequence $\{ S_\tau\}_{\tau}$ converges in probability to a stationary kernel independently on $\tau$, \textit{i.e.}
  \begin{equation}
    S_{\tau}(A_i, A_j) \underset{\mathbb{P}}{\rightarrow} S(A_i, A_j), \; \forall \, \tau > 0, \; \forall \, A_i,A_j.
  \end{equation}
  Hence, the family of NTKs strictly depends on two parameters: activation function and depth of the associated neural network.
  In this paper, with NTK we refer to the infinite-width limit kernel function.

  \item \textbf{Euclidean distance between centroids (IQ)}: $S(A_1, A_2) = \frac{1}{\sqrt{\sigma_1 \| \mathbf{c}_{A_1} - \,  \mathbf{c}_{A_2}\| + \sigma_2}}$, where, even in this case, $ \mathbf{c}_{A_*}$ is the vector of coordinates of the centroid of non-zero elements of discrete observations in $A_*$.
  This kernel function is driven by the physics of the problem at stake, estimating the distance between centroids of the activation region in Euclidean metrics.
\end{itemize}

Each of the kernel functions considered depends on specific hyperparameters that need to be tuned for the given physical application, \textit{e.g.} the variance for RBF kernels, the width and depth of NTK, and the constants $\sigma_1,\sigma_2$ for the IQ kernel. A table summarizing the tuned hyperparameters for each KOL scheme used in the 2D sensitivity analysis (see Section \ref{sec:2Dcase}) is provided in Appendix \ref{app:nomenKOL}.

\section{Numerical results}\label{sec:numerical}
In this section we report the results of the numerical tests performed with FNO and KOL on three different test cases for our problem: a 2D squared domain (Section~\ref{sec:2Dcase}), a 3D slab (Section~\ref{sec:3dcase}) and a realistic ventricle unstructured mesh (Section~\ref{sec:3dunstcase}).

\subsection{Dataset generation details}\label{sec:datasetgen}

\begin{figure}[ht]
    \centering
    \includegraphics[width=\textwidth]{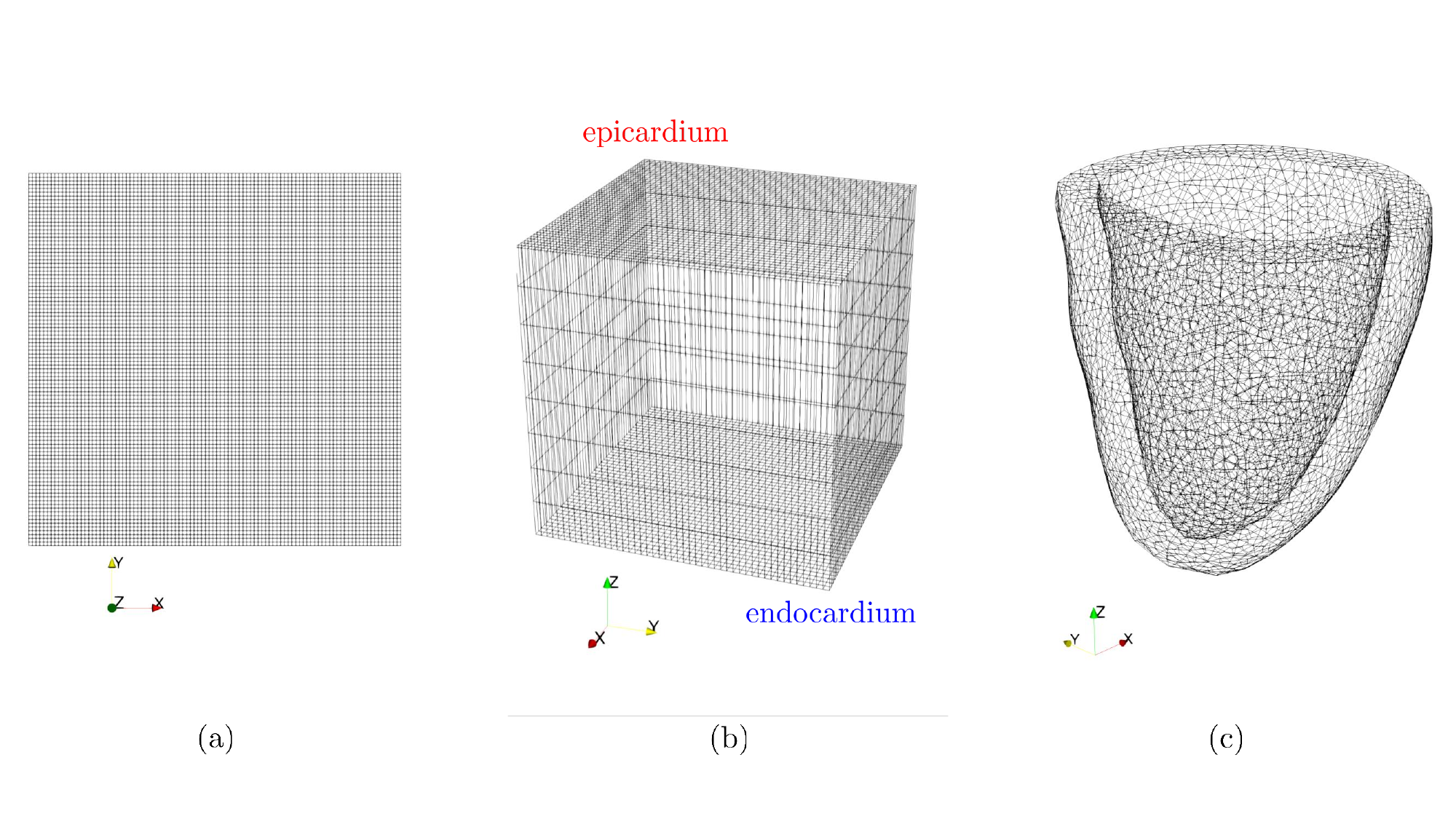}
    \caption{Grids adopted for the numerical simulations: (a) 2D grid, (b) 3D slab and (c) 3D unstructed ventricle.}
    \label{fig:grids}
\end{figure}
To construct the dataset for training operator learning models, we generate input excitations and the corresponding solutions of the 
Monodomain model \eqref{eq:monodomain} discretized by Q1 finite elements, on quadrilateral grids for the 2D case and hexahedral grids for the 3D case (see Figure \ref{fig:grids}). The Monodomain model is coupled with either the Rogers-McCulloch ionic model \cite{rog94} for the 2D case or the Ten Tusscher ionic model \cite{ten06} for the 3D cases, in order to describe the transmembrane potential between the intra- and extracellular domains. Activation and repolarization times are computed as postprocessing of the Monodomain solutions. The applied excitation is defined as a fixed intensity pulse applied for 1ms over a random region of the domain. 
Regarding the diffusion parameters, we have considered $\chi=C_m=1$, and $\mathbf{D}_i = \sigma_1 \mathbf{n}_1 \mathbf{n}_1^T + \sigma_2 \mathbf{n}_2 \mathbf{n}_2^T + \sigma_3 \mathbf{n}_3 \mathbf{n}_3^T$, where the triplet $\{\mathbf{n}_1, \mathbf{n}_2 ,\mathbf{n}_3\}$ represents the fiber orientation.
The conductivities employed in the various scenarios are chosen according to the test cases outlined in \cite{col14}.
For each training we consider an input dataset collecting electrical stimuli (named \texttt{iapps}) and the corresponding output maps of activation times (named \texttt{acti}) and repolarization times (named \texttt{repo}). The name of each dataset is followed by the number of samples it contains, where the repartition between training and testing dataset is 80\%/20\%. The datasets' structure varies depending on the spatial dimension, as detailed below:

\begin{itemize}
    
    \item \textbf{2D case}: we consider $N$ different current stimuli in \texttt{iapps} in the form of $N \times 2$ matrices, collecting the $(x, y)$ coordinates where each pulse is applied. For the basic case, activation and repolarization times are two matrices named of size $N \times n_{\text{no}}^2$, where $n_{\text{no}}$ is the number of discretization points per dimension (100 in our case). In this scenario, we have also considered the case where the cardiac fibers modeled are rotated by 45° counterclockwise, and the corresponding output maps of activation and repolarization times are called \texttt{acti rot} and \texttt{repo rot}, respectively. The sizes of these datasets are the same of the basic case.
    In this case we assume $\mathbf{n}_1 = \mathbf{n}_x$ and $\mathbf{n}_2 = \mathbf{n}_y$, and a counterclockwise rotation of each direction for the rotated case.
    
    \item \textbf{3D slab}: each stimulus is an $N \times n_x \times n_y \times n_z$ binary tensor, where ones indicate pulse application points. We consider $n_x =n_y = 49$ nodes in the $x$ and $y$ directions, $n_z=9$ nodes in the $z$ direction. The activation and repolarization outputs have the same shape of the pulse.
    For the fiber direction we assume $\mathbf{n}_1 = \mathbf{n}_x,\, \mathbf{n}_2 = \mathbf{n}_y$ and $\mathbf{n}_3 = \mathbf{n}_z$.
    
    \item \textbf{3D unstructured ventricle}: the input stimuli are $N \times N \times 3$ matrices, where $N$ is the number of nonzero nodes in the unstructured mesh of the ventricle. Similarly, activation and repolarization outputs are $N \times N_n \times 3$ tensors, with $N_n$ being the total number of unstructured mesh nodes (about 35k degrees of freedom).
    In this case we consider the fiber orientation extracted by the physiological ventricle.
\end{itemize}
These structured and unstructured datasets provide diverse inputs for training operator learning architectures.
In order to evaluate accuracy performance of the trained OL schemes we compute the generalization error as the discrete $L^2$ relative error, namely,
\begin{equation}
  {\mathrm{L^2_{\text{err}}}} = \displaystyle\dfrac{1}{N_p} \sum_{n=1}^{N_p}\displaystyle\dfrac{\| \textbf{y}_p^{n} - \textbf{y}^{n}\|_2}{\| \textbf{y}^{n}\|_2},
  \label{eq:perror}
\end{equation}
where $\textbf{y}^{n}$ is the vector containing the ground truth evaluations of activation/repolarization maps at the different points of the domain whilst $\textbf{y}_p^{n}$ is the corresponding vector of predictions.
For the 3D case representing a 3D ventricle, we remind that the Pearson dissimilarity coefficient $P$ is defined as $P=1-R$, with 
\begin{equation}
R = \frac{\text{Cov}(Y,Y_p)}{\sigma_{Y}\sigma_{Y_p}}.
\end{equation}
Here, $\text{Cov}(Y,Y_p)$ is the covariance between the ground truth of the tested samples and the predictions as flattened vectors (respectively $Y$ and $Y_p$) and $\sigma_{Y,Y_p}$ is the standard deviation relative to the test dataset or the predictions.

\subsection{2D case}\label{sec:2Dcase}
\begin{figure}[H]
    \centering
\begin{minipage}{0.267\textwidth}
    \centering
    \includegraphics[width=\textwidth, trim={0 0 30.4cm 1.5cm}, clip]{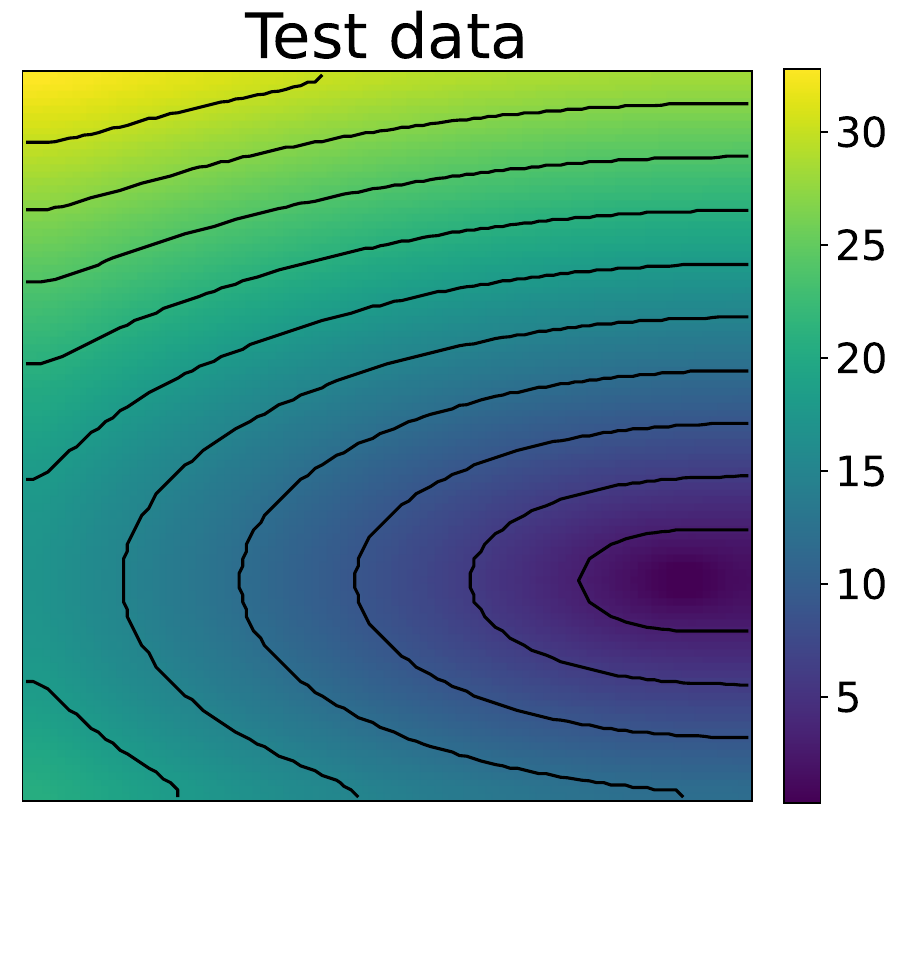}
\end{minipage}
\begin{minipage}{0.534\textwidth}
\hfill\centering
    \includegraphics[width=\textwidth, trim={15.1cm 2cm 0 1.5cm}, clip]{paper_test_2D_001_index_53-eps-converted-to.pdf}
    \caption*{FNO}

\hfill\centering
    \includegraphics[width=\textwidth, trim={15.1cm 2cm 0 1.5cm}, clip]{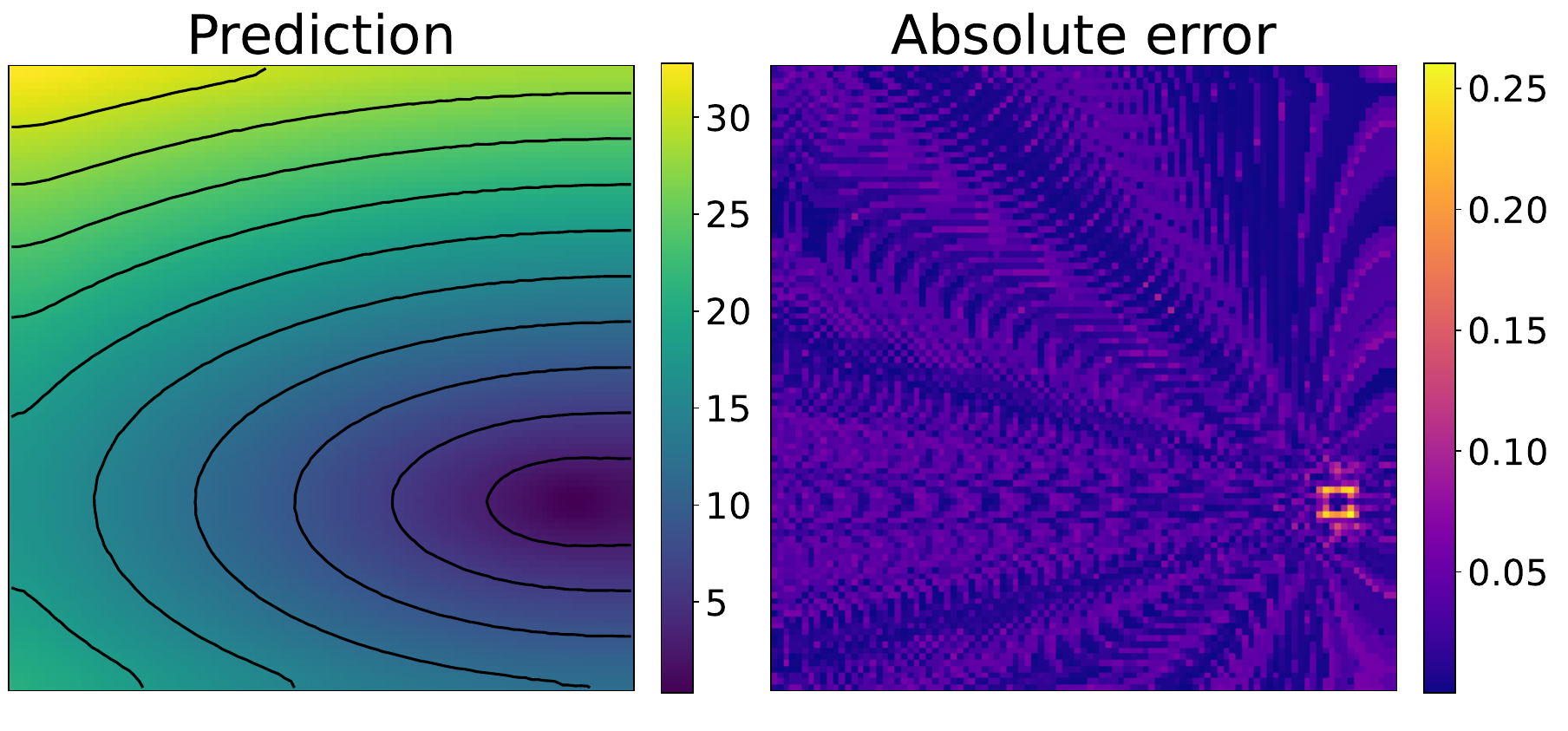}
    \caption*{KOL}
\end{minipage}
\caption{Comparison of FNO (top) and KOL (bottom) activation time predictions for the 2D case (acti with 2000 samples).}
\label{fig:acti2D_comparison}
\end{figure}

\begin{figure}[H]
    \centering
\begin{minipage}{0.267\textwidth}
    \centering
    \includegraphics[width=\textwidth, trim={0 0 30.4cm 1.5cm}, clip]{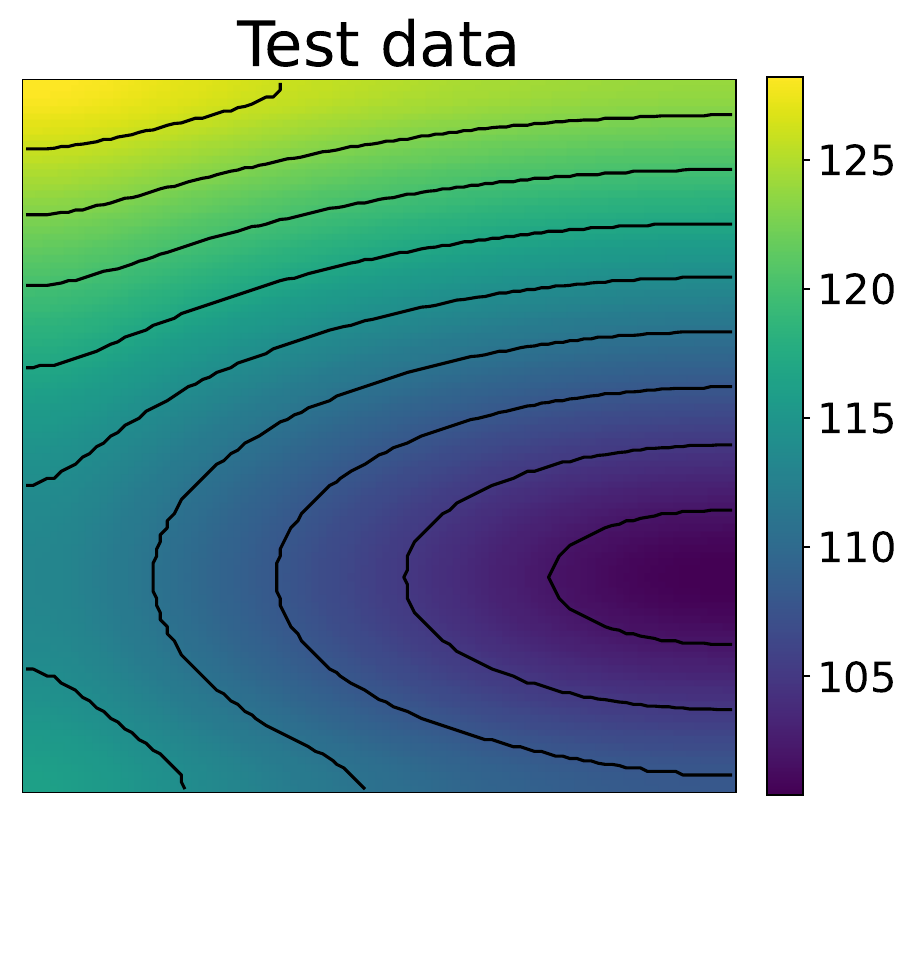}
\end{minipage}
\begin{minipage}{0.534\textwidth}
\hfill\centering
    \includegraphics[width=\textwidth, trim={15.1cm 2cm 0 1.5cm}, clip]{paper_test_2D_007_index_138-eps-converted-to.pdf}
    \caption*{FNO}

\hfill\centering
    \includegraphics[width=\textwidth, trim={15.1cm 2cm 0 1.5cm}, clip]{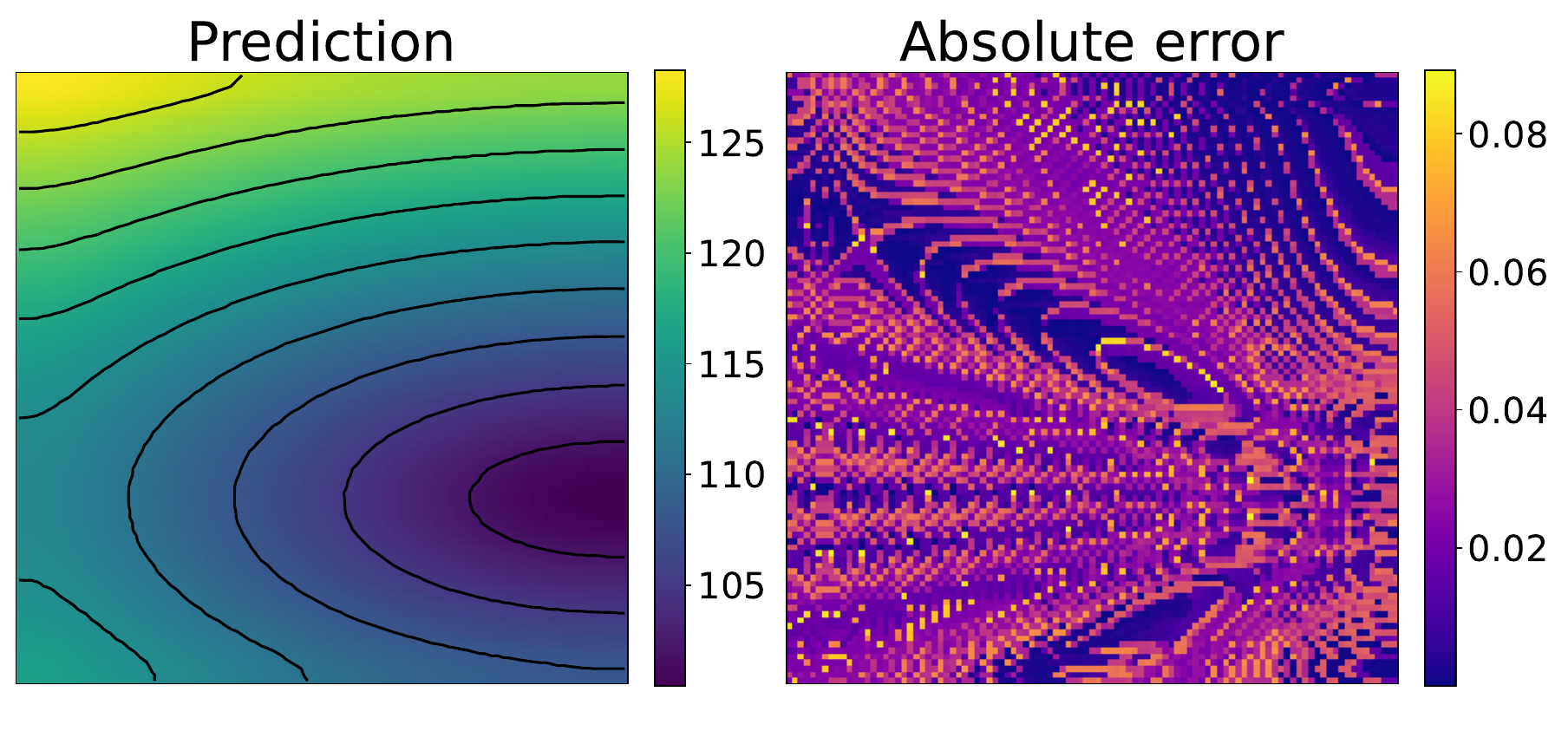}
    \caption*{KOL}
\end{minipage}
\caption{Comparison of FNO (top) and KOL (bottom) repolarization time predictions for the 2D case (repo with 2000 samples).}
\label{fig:repo2D_comparison}
\end{figure}

The experiments conducted on the 2D domain $\Omega = (0,1)^2$ using both architectures provide insightful observations on the model's performance across different datasets, hyperparameters setup, and fiber configurations.
Figures \ref{fig:acti2D_comparison} and \ref{fig:repo2D_comparison} depict test samples reconstructed with FNO (top) and KOL (bottom) for activation and repolarization times.

The FNO architecture employed comprises four Fourier layers with 16 Fourier modes along the $x$-axis and 4 modes along the $y$-axis, with a total number of approximately 532k trainable parameters.
In Tables~\ref{table:2D_results} and \ref{table:2D_results_FNO_appendix}, results indicate that the use of a tailored learning rate reduction policy (\texttt{reduce-} \texttt{OnPlateau}), where the learning rate is reduced by a factor of 0.95 whenever the test loss is not decreasing, consistently outperforms training without such a policy. 
For instance, in the activation dataset \texttt{acti} with 3000 samples, the test error is reduced from $2.82 \times 10^{-3}$ (no policy) to $2.66 \times 10^{-3}$ when the policy is applied.
In the FNO results the uncertainty bands arise from considering different trained architectures with various Kaiming normal initializations \cite{he15} of the trainable parameters. 

The rotated fiber tests, where a 45° rotation of the myocytes is applied on the domain, generally show larger test errors compared to their non-rotated counterparts with the same number of training samples (cfr. Table \ref{table:2D_results_FNO_appendix} and \ref{table:2D_results_kol_acti_appendix}). This observation suggests that the FNO architecture might be sensitive to the structural alignment of input features.
Moreover, we notice that increasing the dataset size improves the performance, as evident in both activation and repolarization datasets. For activation tests with 2000 and 3000 samples, the test error drops, highlighting the benefit of having more training data for learning the spatial features of the solution. In terms of computational efficiency, for FNO the GPU memory consumption remained stable at approximately 5.6GB for 2000-sample datasets and increased slightly to 5.79GB for 3000-sample datasets. Training times scaled proportionally with dataset size, from 48 minutes for 2000 samples to 72 minutes for 3000 samples.
We notice a general low absolute error, with some areas for the activation case propagating orthogonally with respect to level curves. This error distribution was observed also for the KOL case in Figure~\ref{fig:acti2D_comparison} (bottom, right), but a theoretical explanation of such a pattern remains unknown. A more uniform pattern is observed instead for the repolarization case.

We conducted a large number of numerical simulations using the same training and testing datasets employed for FNO architectures to evaluate the performance of KOL (cfr. Tables \ref{table:2D_results}, \ref{table:2D_results_kol_acti_appendix}, and \ref{table:2D_results_kol_repo_appendix}). The primary objective of this sensitivity analysis was to assess the impact of kernel selection, a critical factor that is highly problem-dependent and can significantly affect prediction quality.

From Tables \ref{table:2D_results_kol_acti_appendix} and \ref{table:2D_results_kol_repo_appendix}, we observe that IQ kernel-based strategies yield generalization errors at least two orders of magnitude lower than those using NTK or RBF. This improved performance is attributed to the IQ kernel’s efficiency in computing correlations of compact support indicator functions, which accurately represent activation regions. Notably, this advantage persists even when increasing the training set size for both activation and repolarization reconstructions.
Additionally, KOL exhibits sensitivity to structural alignments of input features. Specifically, tests with rotated fibers achieve higher accuracy than unrotated counterparts. We observe that KOL significantly reduces training times with respect to FNO, decreasing from thousands of seconds to just few hundreds, as it requires solving a single symmetric, positive definite linear system rather than an iterative optimization process. However, as training size increases, the system's condition number rises which may impact on testing performances. Therefore, in this case, the use of tailored preconditioning strategies is crucial (cfr. \cite{mea20, shi24}).
Regarding computational efficiency, CPU memory consumption scales with training size but remains close to 1 GB. It is important to note that the accuracy and training time improvements offered by KOL must be considered in light of the non-negligible time required for kernel selection.
Furthermore, we note that KOL, equipped with the chosen deterministic kernels, produces fully deterministic predictions: therefore, unlike FNO, it does not exhibit uncertainty bands due to hyperparameter initialization.

\begin{table}
    \centering
    \small
    \begin{tblr}{
        colspec = {ccccc}
    }
        \toprule
        \SetCell[r=4]{c} {FNO \\ with lr policy \\ \texttt{reduceOnPlateau}} 
        & Test Error & Dataset (size) & {GPU \\ Memory} & {Training \\ Time} \\
        \midrule
        & 2.66E-03 $\pm$ 2.27E-04 & \texttt{acti} (3000) & 5.79GB & 72 min \\
        & 3.33E-03 $\pm$ 3.13E-04 & \texttt{acti rot} (2000) & 5.62GB & 48 min \\
        & 3.13E-03 $\pm$ 2.33E-04 & \texttt{repo} (3000) & 5.79GB & 72 min \\
        & 3.53E-03 $\pm$ 3.63E-04 & \texttt{repo rot} (2000) & 5.62GB & 48 min \\
        \midrule
        \midrule
        \SetCell[r=4]{c} {KOL \\ with \texttt{iq4} kernel}
        & Test Error & Dataset (size) & {CPU \\ Memory}& {Training\\Time} \\
        \midrule
         & 9.52E-04 & \texttt{acti} (3000) & 1.15GB & 10 min \\
        & 9.34E-04 & \texttt{acti rot} (2000)& 0.91GB & 8 min \\
        & 4.74E-04 & \texttt{repo} (3000) & 1.15GB & 10 min \\
        & 4.69E-04 & \texttt{repo rot} (2000)& 0.91GB & 6 min \\
        \bottomrule
    \end{tblr}
    \caption{Performance comparison of FNO (reduce-on-plateau learning rate policy) and KOL (\texttt{iq4} kernel) methods on 2D datasets.}
    \label{table:2D_results}
\end{table}


\subsection{3D slab}\label{sec:3dcase}
\begin{figure}
    \centering
    \includegraphics[width=0.45\linewidth]{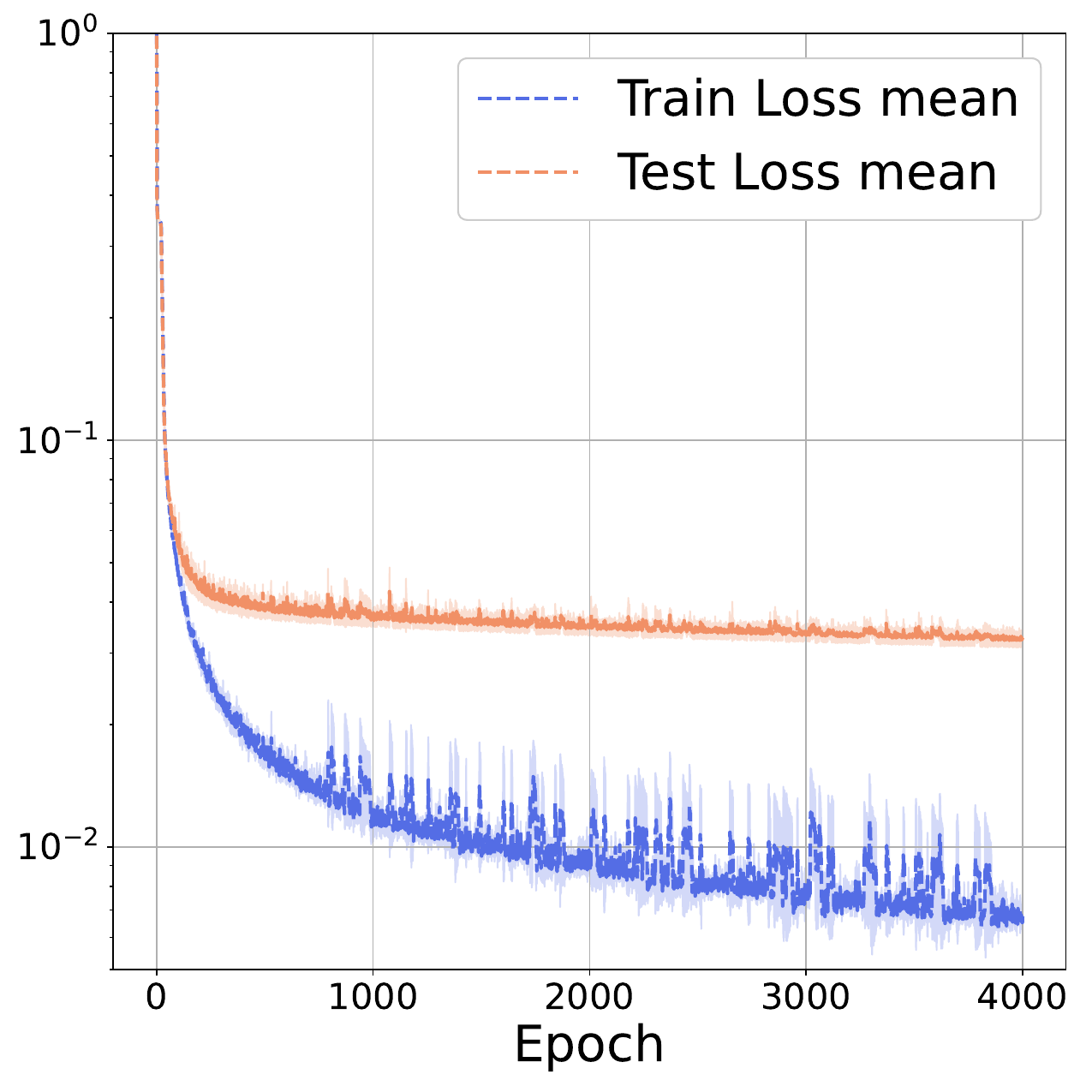}
    \caption{FNO loss plot for the 3D case (\texttt{acti} 2000 dataset). Mean train and test loss of three different randomly initialized models (dashed line). Standard deviation over the three models for each epoch is also reported (light shadow).}
    \label{fig:fno_loss}
\end{figure}

\begin{figure}
\centering
\begin{subfigure}{0.24\textwidth}
    \centering
    \includegraphics[width=\textwidth]{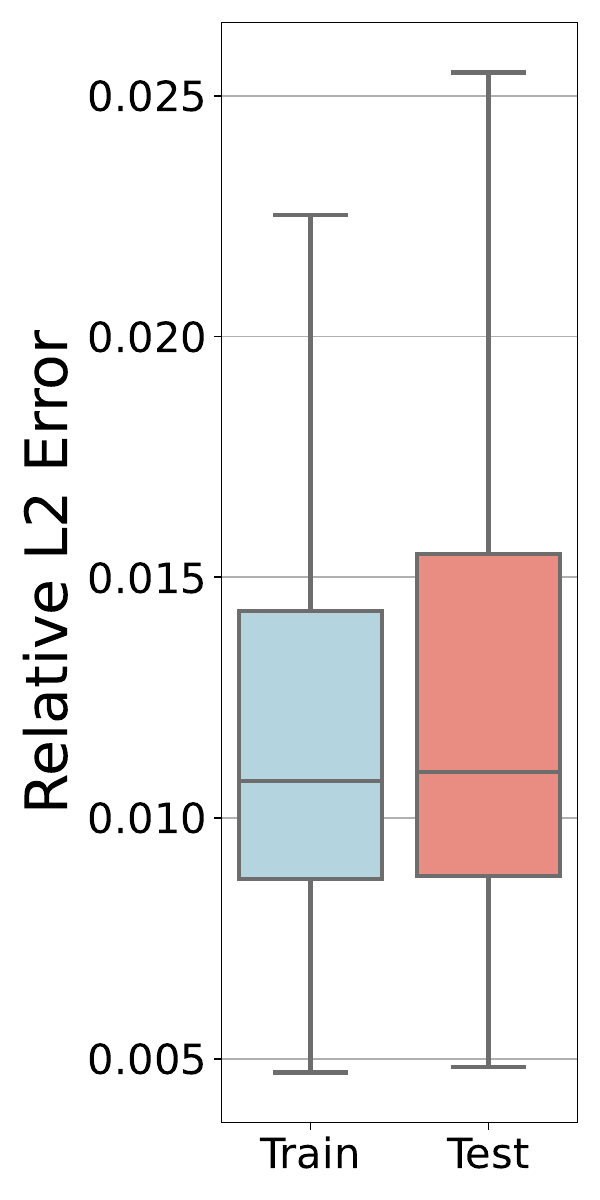}
    \caption{FNO box plot.}
    \label{fig:fno_box}
\end{subfigure}
    \begin{subfigure}{0.55\textwidth}
    \centering
    \includegraphics[width=\textwidth]{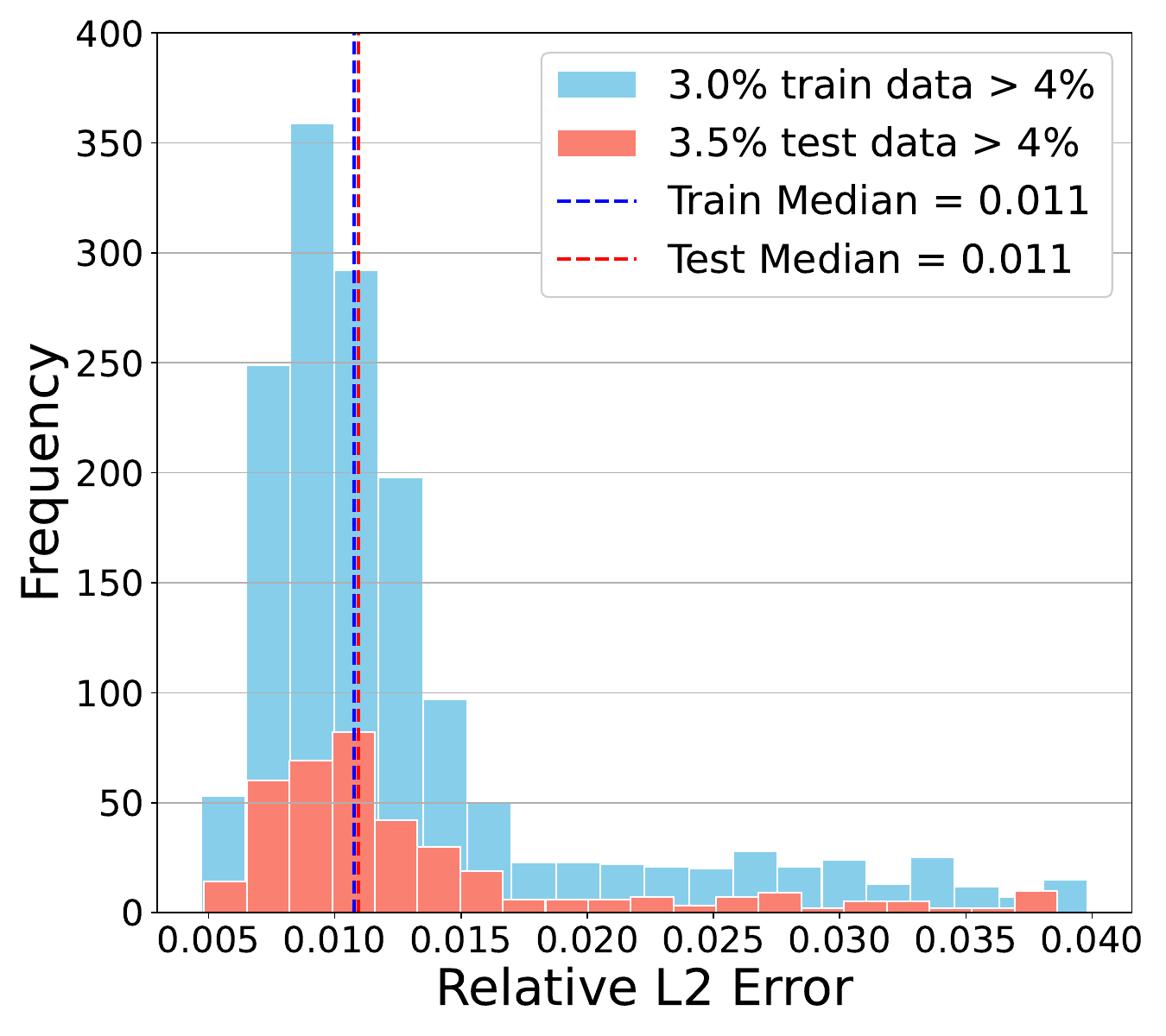}
    \caption{FNO histogram.}
    \label{fig:fno_hist}
\end{subfigure}
\caption{FNO box plot and histogram for the 3D dataset \texttt{acti2000} relative to the best model trained (outliers not shown). For the training set, 3\% of the data have a relative L2 error greater than 4\%, while for the test set, 3.5\% of the data exceed this threshold.}
\label{fig:statsFNO}
\end{figure}

\begin{figure}
\centering
\begin{subfigure}{0.18\textwidth}
    \centering
    \includegraphics[width=\textwidth]{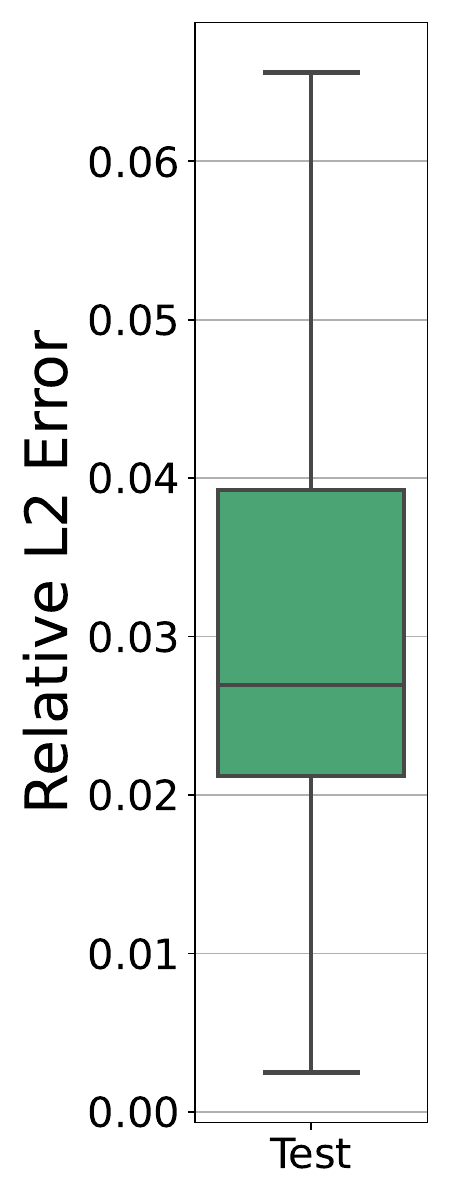}
    \caption{KOL box plot.}
    \label{fig:kol_box}
\end{subfigure}
    \begin{subfigure}{0.55\textwidth}
    \centering
    \includegraphics[width=\textwidth]{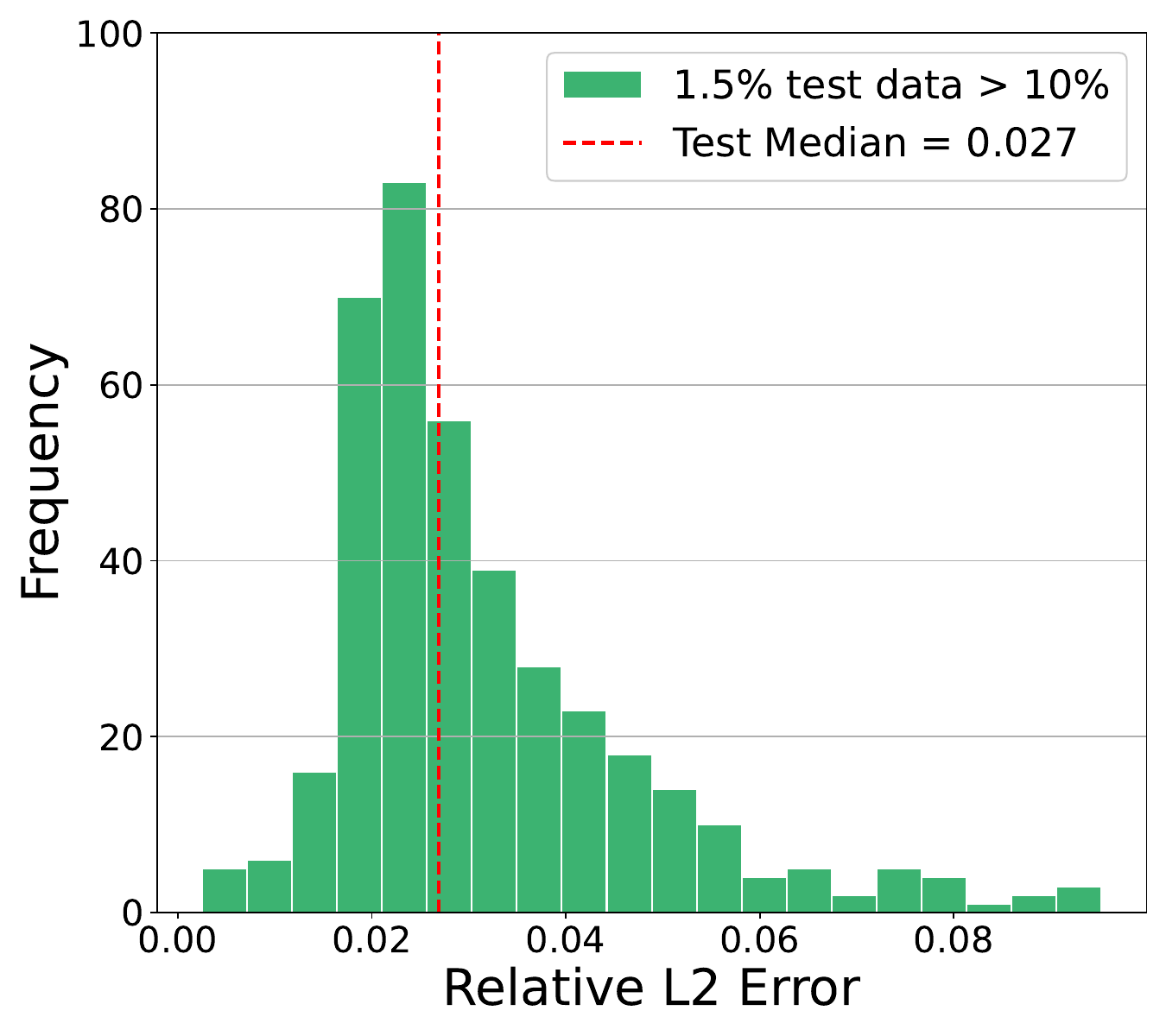}
    \caption{KOL histogram.}
    \label{fig:kol_hist}
\end{subfigure}
\caption{KOL box plot and histogram for 3D dataset \texttt{acti} 2000 relative to the best model trained (outliers not shown). Training results are not shown since we achieve machine precision. For the test set, 1.5\% of the data have a relative L2 error greater than 10\%.}
\label{fig:statsKOL}
\end{figure}

Additional experiments conducted on a 3D slab in $(0,1)^3$ have been conducted using the two operator learning approaches discussed. As will be shown, the performance of FNO and KOL remained robust despite the increased dimensionality of this test case.
The architecture employed for FNO comprises four Fourier layers with 16 Fourier modes along the $x$-axis, 8 along the $y$-axis and 4 along the $z$-axis, resulting in approximately 8.4 million trainable parameters. Given the improved performance observed for the 2D case, all experiments used the \texttt{reduceOnPlateau} learning rate policy, where the learning rate decreases upon stagnation of the validation loss.
Instead, we endow KOL with \texttt{iq4} kernel (cfr. Table \ref{table:nomenKOL_appendix}) following the sensitivity analysis of the 2D case.
In Table~\ref{table:3D_results}, the results indicate that increasing the dataset size significantly enhances prediction accuracy for both models. For the activation dataset, the FNO test error decreased from $4.15 \times 10^{-2}$ to $3.27 \times 10^{-2}$ when the sample size increased from 1000 to 2000. A similar trend is observed for repolarization, where the test error dropped from $8.46 \times 10^{-3}$ to $6.91 \times 10^{-3}$.
Conversely, KOL exhibits the opposite behavior: as the dataset size grows, the conditioning number of the SDP system increases, and, therefore testing error increases (\textit{e.g.} from $ 5.33 \times 10^{-3}$ at 1000 dataset size to $ 6.19 \times 10^{-3}$ at 2000 dataset size for repolarization).
Nonetheless, test errors for KOL remain consistently lower than those of FNO in both activation and repolarization cases.
However, as reported in Figure~\ref{fig:statsFNO}, most of the test data predictions for the \texttt{acti 2000} dataset are distributed under 1\% of the relative L2 error, despite having around 3.5\% of the data as outliers for the FNO case.
Similar results are obtained for KOL in Figure~\ref{fig:statsKOL}. For the \texttt{repo 2000} dataset the reader can refer to Figures~\ref{fig:statsFNOrepo}-\ref{fig:statsKOLrepo}. In Figure~\ref{fig:fno_loss} we have also reported the training and test loss decaying plot for FNO. 

\begin{figure}
    \centering
    \includegraphics[width=0.75\textwidth]{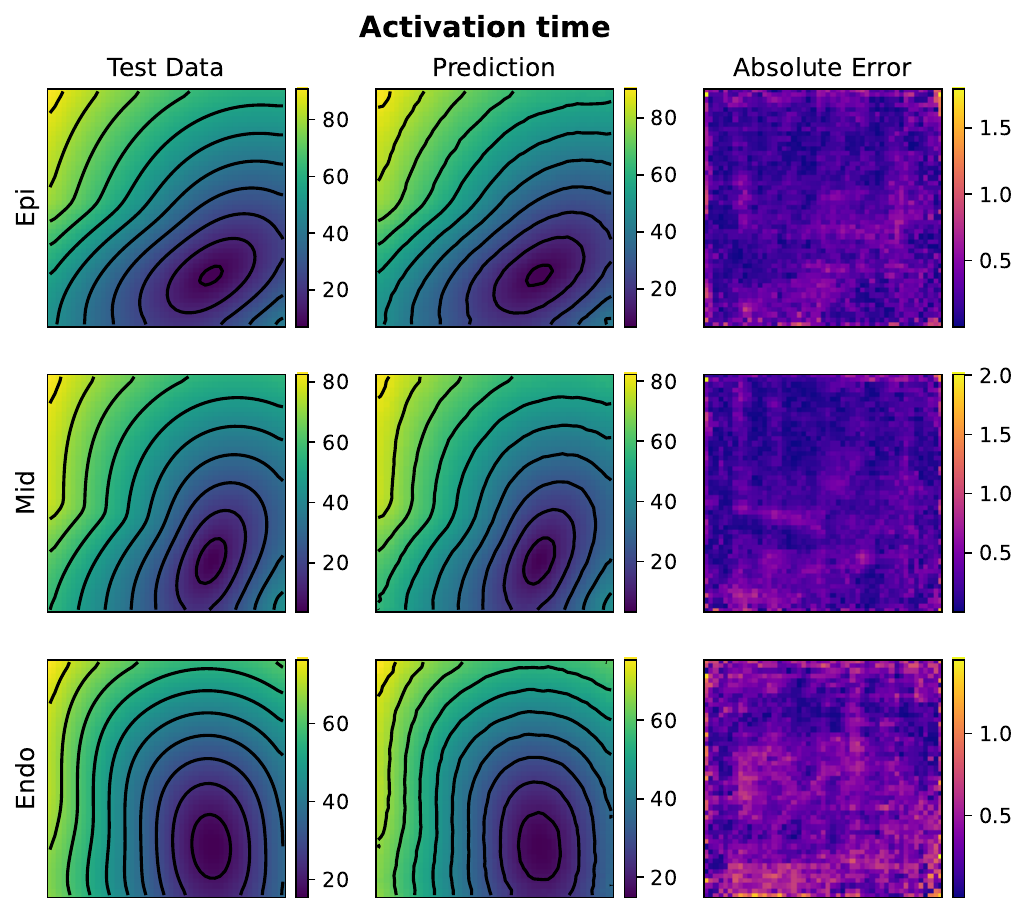}
    \caption{Example of FNO activation time prediction for the 3D case (acti with 2000 samples). The picture represents three slices of tissue: epicardium (top), middle (center) and endocardium (bottom).}
    \label{fig:acti3D}
\end{figure}

\begin{figure}
    \centering
    \includegraphics[width=0.75\textwidth]{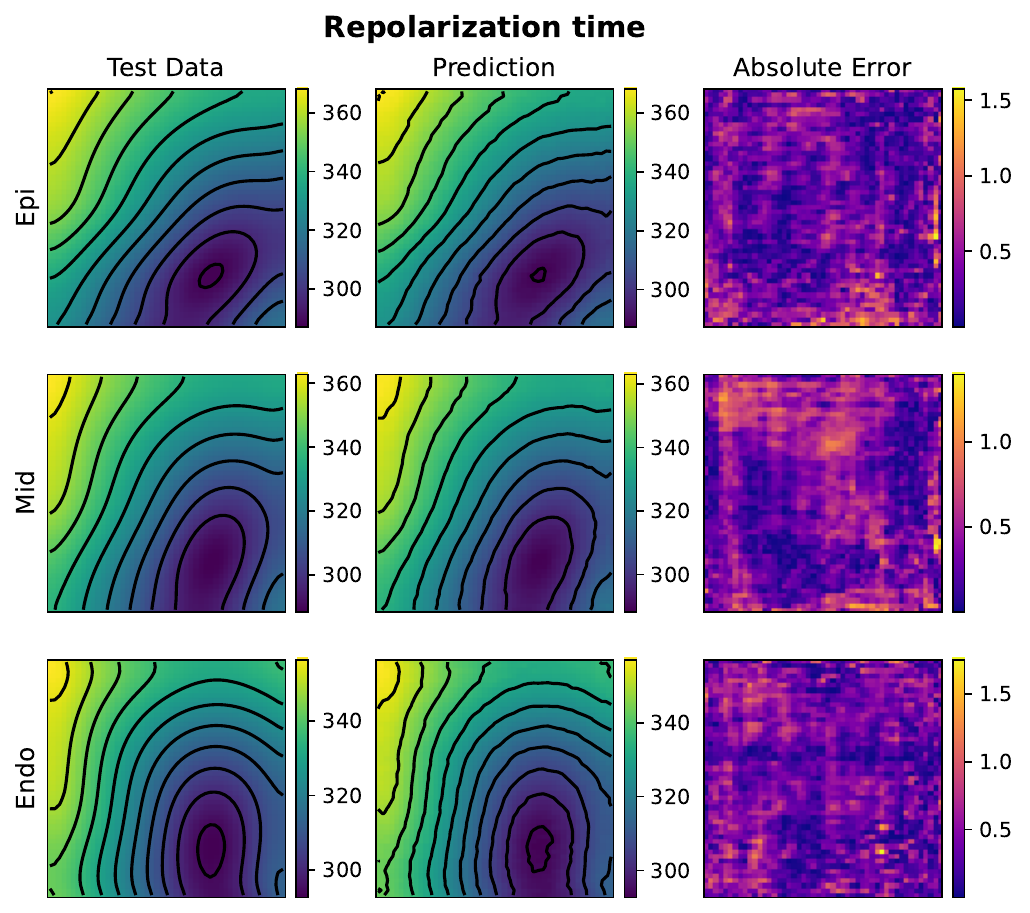}
    \caption{Example of FNO repolarization time prediction for the 3D case (repo with 2000 samples). The picture represents three slices of tissue: epicardium (top), middle (center) and endocardium (bottom).}
    \label{fig:repo3D}
\end{figure}

\begin{figure}
    \centering
    \includegraphics[width=0.75\textwidth]{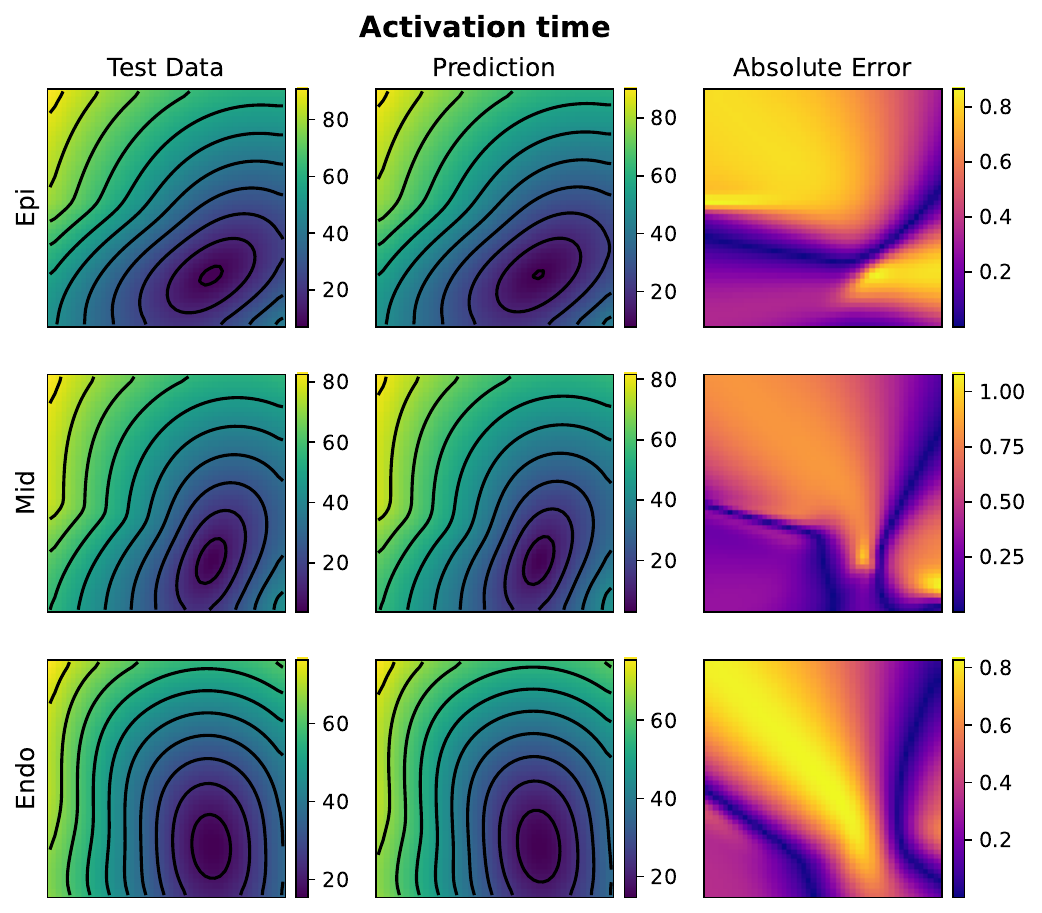}
    \caption{Example of KOL activation time prediction for the 3D case (acti with 2000 samples). The picture represents three slices of tissue: epicardium (top), middle (center) and endocardium (bottom).}
    \label{fig:acti3D_kol}
\end{figure}

\begin{figure}
    \centering
    \includegraphics[width=0.75\textwidth]{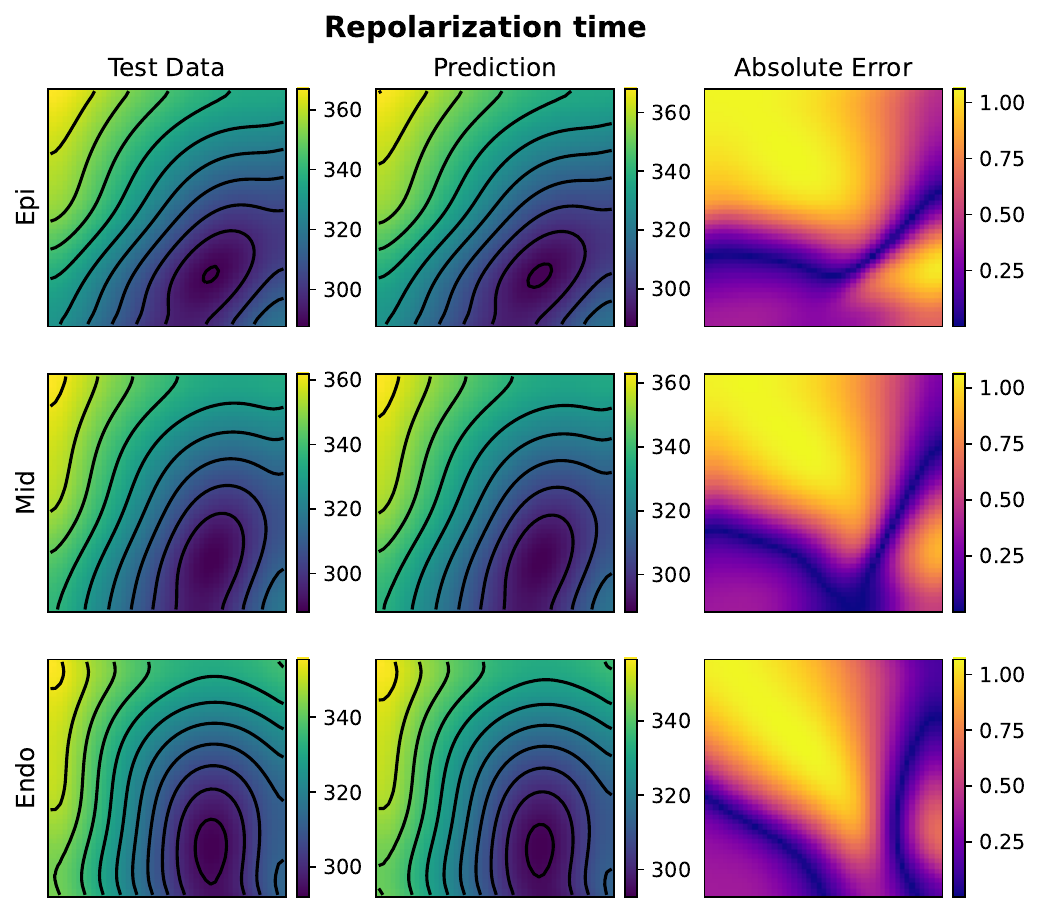}
    \caption{Example of KOL repolarization time prediction for the 3D case (repo with 2000 samples). The picture represents three slices of tissue: epicardium (top), middle (center) and endocardium (bottom).}
    \label{fig:repo3D_kol}
\end{figure}

Compared to 2D tests, computational costs rise substantially in 3D. For FNO, GPU memory consumption increased from 13.66GB (1000 samples) to 14.09GB (2000 samples), whereas KOL required 1.67GB and 2.56GB, respectively. Training times scaled linearly for both models, with FNO requiring 94 minutes for 1000 samples and 188 minutes for 2000, while KOL completed training in 211 seconds and 427 seconds for the same dataset sizes.
Figures~\ref{fig:acti3D}, \ref{fig:repo3D}, \ref{fig:acti3D_kol} and \ref{fig:repo3D_kol} illustrate activation and repolarization time predictions for both models, along with the corresponding high fidelity solutions and absolute errors across three slices of the slab domain representing the endocardium, epicardium and the intermediate slice. The applied stimulus belongs to the epicardium surface for the reconstruction of activation times, whilst it belongs to an intermediate sheet between the endocardium and the middle of the slab for the repolarization case.
In both cases, prediction error increases as far as the distance from applied stimuli increase.
While FNO exhibits a nearly uniform absolute error distribution slightly higher than in the 2D case, KOL’s errors tend to be concentrated near the activation region and propagate orthogonally to the level curves of the activation (or repolarization) times.

\begin{table}
    \centering
    \small
    \begin{tblr}{
        colspec = {ccccc}
    }
        \toprule
        \SetCell[r=4]{c} {FNO \\ with lr policy \\ \texttt{reduceOnPlateau}} 
        & Test Error & Dataset (size) & {GPU \\ Memory} & {Training \\ Time} \\
        \midrule
        & 4.15E-02 $\pm$ 2.02E-03 & \texttt{acti} (1000) & 13.66GB & 94 min \\
        & 3.27E-02 $\pm$ 5.86E-04 & \texttt{acti} (2000) & 14.09GB & 188 min \\
        & 8.46E-03 $\pm$ 6.51E-04 & \texttt{repo} (1000) & 13.66GB & 94 min \\
        & 6.91E-03 $\pm$ 3.44E-04 & \texttt{repo} (2000) & 14.09GB & 188 min \\
        \midrule
        \midrule
        \SetCell[r=4]{c} {KOL \\ with \texttt{iq4} kernel} 
        & Test Error & Dataset (size) & {CPU \\ Memory}& {Training\\Time} \\
        \midrule
        & 1.67E-02 & \texttt{acti} (1000) & 1.67GB & 3 min \\
        & 1.82E-02 & \texttt{acti} (2000) & 2.56GB & 7 min \\
        & 5.33E-03 & \texttt{repo} (1000) & 1.67GB & 3 min \\
        & 6.19E-03 & \texttt{repo} (2000) & 2.56GB & 7 min \\
        \bottomrule
    \end{tblr}
    \caption{Performance comparison of FNO and KOL methods on 3D datasets.}
    \label{table:3D_results}
\end{table}

\subsection{3D unstructured ventricle}\label{sec:3dunstcase}
In the last test case considered, we solved our problem on an unstructured mesh representing a human ventricle, which consists of approximately 35k degrees of freedom.
In this case, the FNO was adapted to handle irregular domains. In particular, we implemented a Non-Uniform Discrete Fourier Transform (NUDFT) following the approach presented in~\cite{lin23}.
On the other hand, KOL did not require any specific modifications to operate in this unstructured setting.
Since the input is defined on an unstructured grid, we modified its representation accordingly. 
Specifically, each dataset consists of $N$ samples, where each sample corresponds to a vector of size $N_{\text{nodes}}$ with binary values: a value of 1 at a given node indicates the presence of an external stimulus, while 0 represents the absence of stimulation. Given that all stimuli have the same pulse intensity and duration, the input applied current was reformulated as a matrix of dimension $N \times N_{\text{stim}} \times 3$, where $N_{\text{stim}}$ represents the maximum number of stimulated nodes across all samples, with padding applied where necessary. The numerical solutions for the activation and the repolarization times were structured as tensors of size $N \times N_{\text{nodes}} \times 3$, capturing the time evolution of the solution on the unstructured mesh.

The results in Table~\ref{table:3D_unst_results} provide insights into the performance of FNO and KOL on this complex domain. The FNO architecture employed consisted of Fourier layers with four Fourier modes in each spatial direction ($x$, $y$, and $z$), leading to trainable parameters ranging between 10.8M and 11.3M. We applied a \texttt{reduceOnPlateau} learning rate policy, as it consistently yielded better predictions in previous tests by dynamically adjusting the learning rate when the test loss plateaued.
The architecture's depth and width also influenced performance. Notably, for $N = 2000$, increasing the model width from 2 to 16 substantially improved accuracy, reducing the test error from $1.44 \times 10^{-1}$ to $7.18 \times 10^{-2}$. However, beyond a certain point, further increases in width did not yield significant improvements, with test errors fluctuating for widths above 32.
Similarly, deeper architectures ($L > 1$) did not always lead to better performance, indicating that an optimal hyperparameters selection is crucial for balancing expressivity and generalization.
In Table \ref{table:3D_results} we reported the most accurate performances of FNOs obtained considering $L=3$.

\begin{figure}[H]
    \centering
\begin{subfigure}{\textwidth}
    \centering                    
    \includegraphics[width=0.25\textwidth, trim={3cm 2cm 50cm 9cm}, clip]{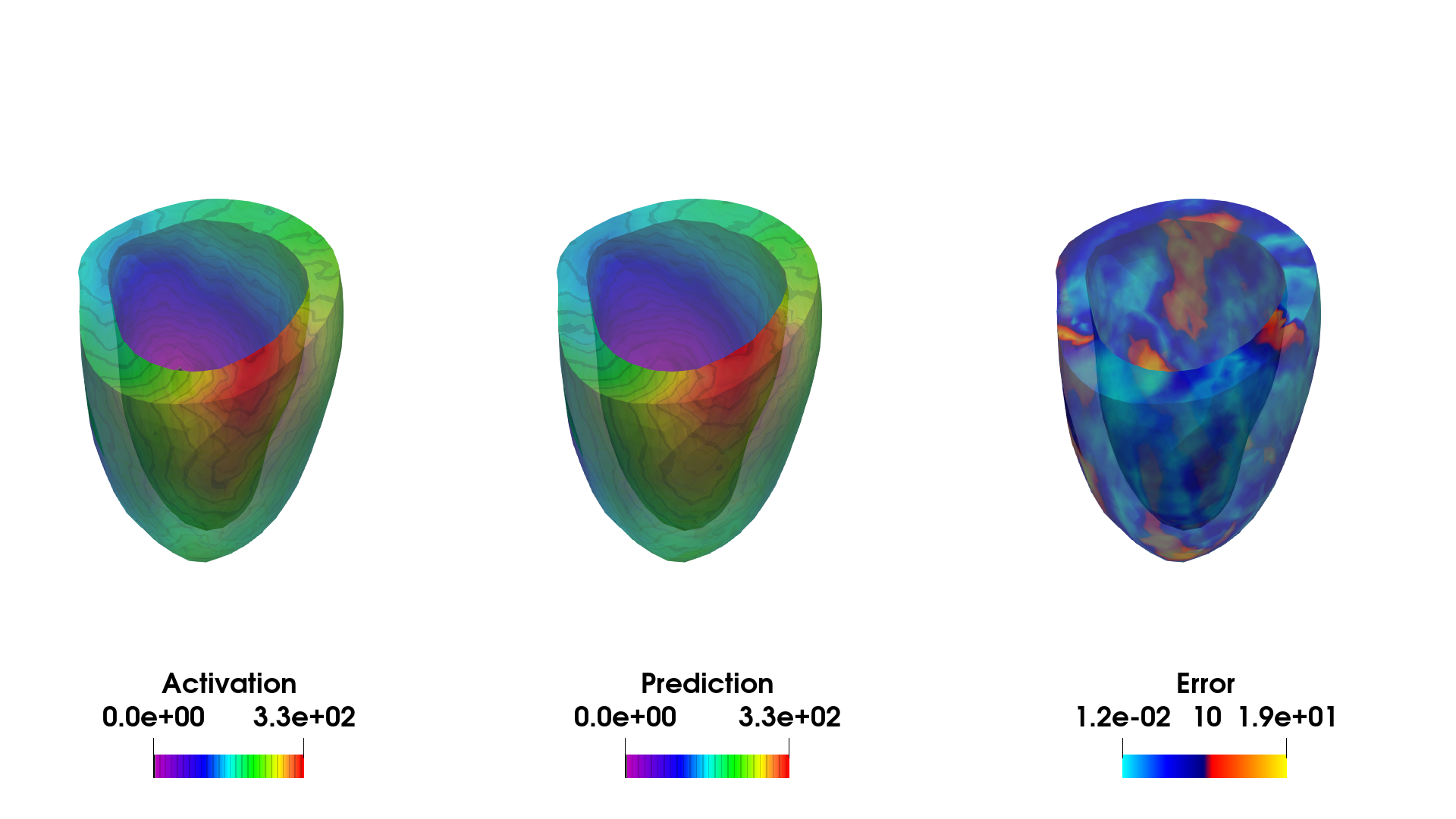}%
    \includegraphics[width=0.25\textwidth, trim={25cm 2cm 28cm 9cm}, clip]{acti_fno_3dunst.png}%
    \includegraphics[width=0.25\textwidth, trim={48cm 2cm 5cm 9cm}, clip]{acti_fno_3dunst.png}
    \caption{Activation}
    \label{fig:acti3D_un}
\end{subfigure}
    \begin{subfigure}{\textwidth}
    \centering
    \includegraphics[width=0.25\textwidth, trim={3cm 2cm 50cm 9cm}, clip]{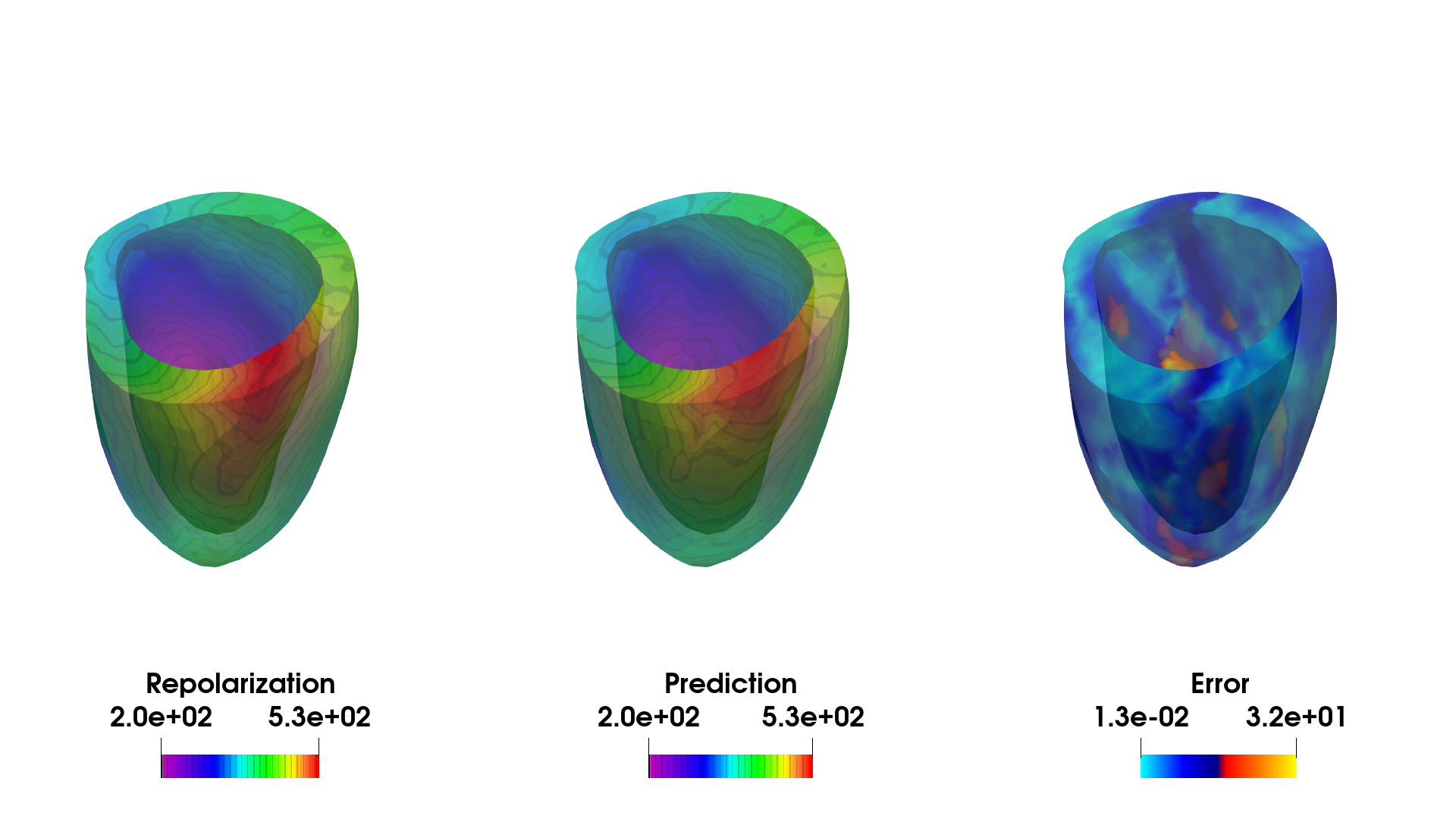}%
    \includegraphics[width=0.25\textwidth, trim={25.5cm 2cm 27.5cm 9cm}, clip]{repo_fno_3dunst.png}%
    \includegraphics[width=0.25\textwidth, trim={48cm 2cm 5cm 9cm}, clip]{repo_fno_3dunst.png}
    \caption{Repolarization}
    \label{fig:repo3D_un}
    \end{subfigure}
    \caption{Example of FNO predictions for the 3D unstructured case: (a) Activation times (acti 2000), (b) Repolarization times (repo 2000).}
    \label{fig:3D_un_FNO}
\end{figure}

\begin{figure}[H]
    \centering
\begin{subfigure}{\textwidth}
    \centering
    \includegraphics[width=0.25\textwidth, trim={3cm 2cm 50cm 8cm}, clip]{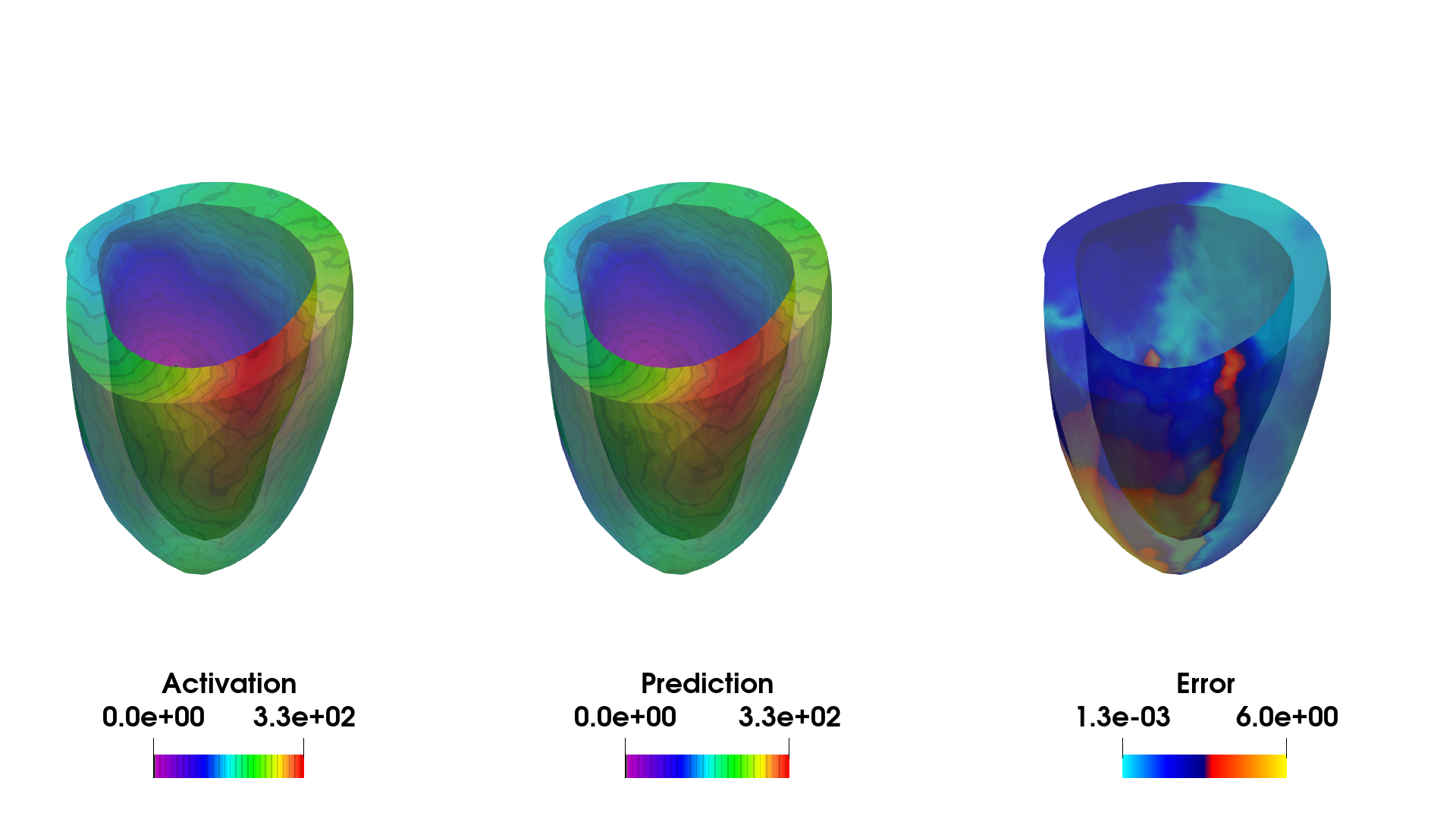}%
    \includegraphics[width=0.25\textwidth, trim={25cm 2cm 28cm 8cm}, clip]{acti_kol_3dunst.png}%
    \includegraphics[width=0.25\textwidth, trim={48cm 2cm 5cm 8cm}, clip]{acti_kol_3dunst.png}
    \caption{Activation}
    \label{fig:acti3D_un_kol}
\end{subfigure}
    \begin{subfigure}{\textwidth}
    \centering
    \includegraphics[width=0.25\textwidth, trim={3cm 2cm 50cm 9cm}, clip]{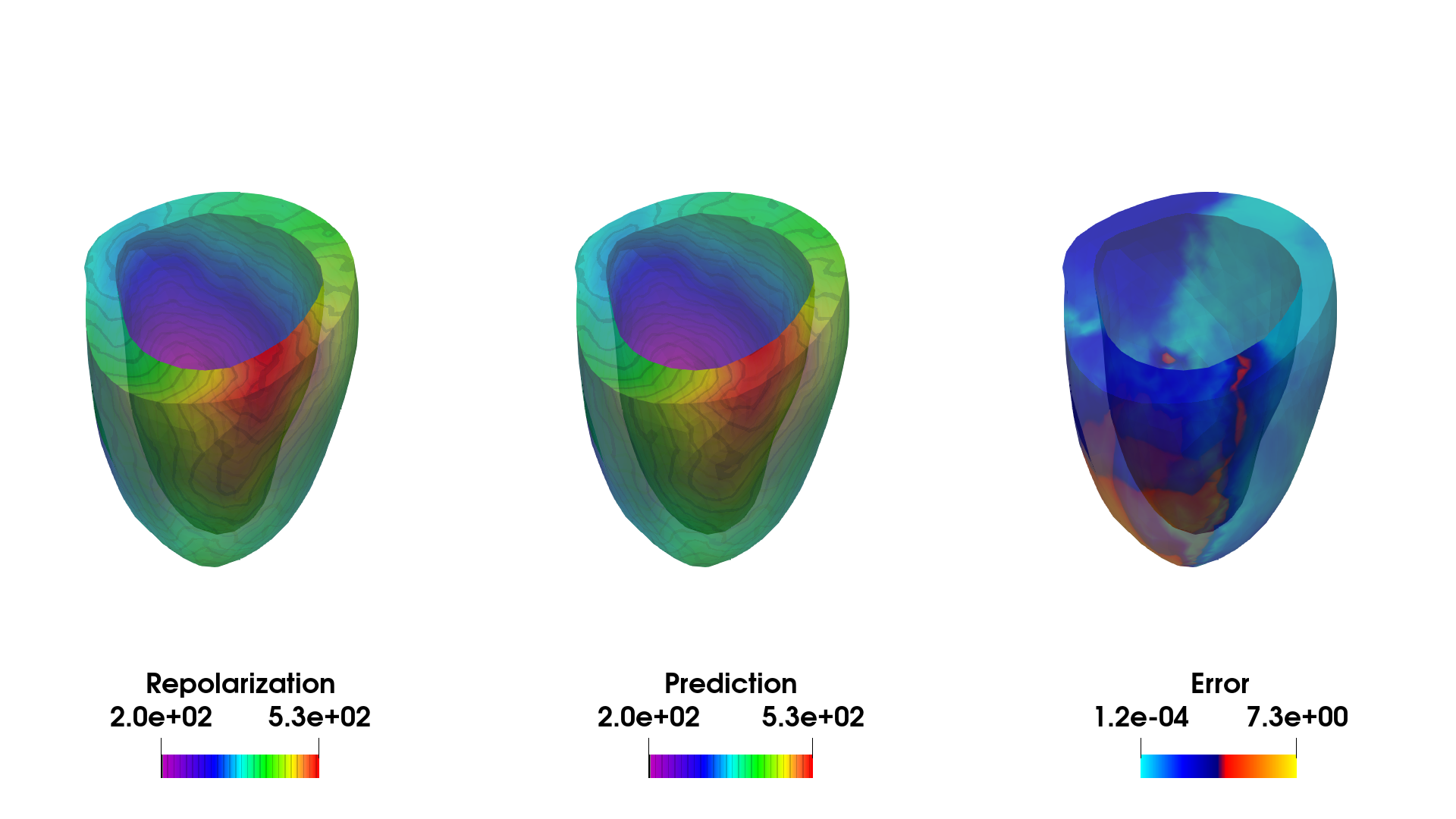}%
    \includegraphics[width=0.25\textwidth, trim={25.5cm 2cm 27.5cm 9cm}, clip]{repo_kol_3dunst.png}%
    \includegraphics[width=0.25\textwidth, trim={48cm 2cm 5cm 9cm}, clip]{repo_kol_3dunst.png}
    \caption{Repolarization}
    \label{fig:repo3D_un_kol}
    \end{subfigure}
    \caption{Example of KOL predictions for the 3D unstructured case: (a) Activation times (acti 2000), (b) Repolarization times (repo 2000).}
    \label{fig:3D_un_KOL}
\end{figure}

Also in this case, KOL (endowed with \texttt{iq4} kernel) significantly outperforms FNO in terms of test error, achieving a test error as low as $5.29 \times 10^{-3}$ on the repolarization dataset for $N = 2000$, compared to FNO's $4.86 \times 10^{-2}$. The difference is even more pronounced for the activation dataset, where KOL attains a test error of $1.22 \times 10^{-2}$, while FNO's best result remains at $7.15 \times 10^{-2}$. Additionally, KOL exhibits lower Pearson dissimilarity values, indicating a better linear correlation between the ground truth test data and the corresponding predictions.
Solutions for a specific test case, for both activation and repolarization, are shown in Figures~\ref{fig:3D_un_FNO} and \ref{fig:3D_un_KOL}. The absolute error plots highlight the lower error achieved by KOL, with a maximum absolute error of 4 ms, compared to larger regions reaching approximately 8.6 ms in the FNO case. However, both architectures successfully captured the qualitative distribution of activation and repolarization times.

Computationally, FNO has an advantage in terms of inference time. A single prediction with FNO on the larger dataset, consisting of 2000 samples with 1600 for training and 400 for testing, takes approximately $8.5$ ms, whereas KOL requires $56$ ms for the same task.
Both methods outperform the solution of a single Monodomain model implemented using the PETSc library \cite{bal19} on a node equipped with an Intel(R) Xeon(R) Silver 4316 CPU @ 2.30GHz and two Nvidia H100 GPUs, which requires about 4 minutes on 2 cores with GPU acceleration. In order to have a fair comparison from a user's point of view, we have performed the single prediction timing tests on an laptop equipped with an Apple M1 Pro chip, while the trainings have been performed on a machine equipped with an NVIDIA Quadro RTX 5000 GPU for FNO, while on an Intel i7 machine for KOL.

Furthermore, in this case training KOL is significantly faster than FNO, requiring only a few minutes compared to FNO's training time, which ranges from 58 to 84 minutes depending on the configuration. KOL is also much more memory efficient, using just over 1GB of CPU memory, compared to FNO's GPU memory consumption of around 6.1GB. These results suggest that KOL is lighter, faster and more accurate than FNO, even though we expect the latter to be more time-efficient when predicting a large number of occurrences is required.

\begin{table}
    \centering
    \small
    \begin{tblr}{
        colspec = {ccccccc}
    }
        \toprule
        \SetCell[r=3]{c} {FNO \\ with lr policy \\ \texttt{reduceOnPlateau}} 
        & Test Error & Dataset (size) & {GPU \\ Memory} & {Training \\ Time} & {Testing \\ Time} & {Pearson \\ (test)} \\
        \midrule
        & 7.15E-02 $\pm$ 9.84E-03 & \texttt{acti} (2000)   & 6.13GB & 80 min & 8.5E-03 sec & 4.3E-03 \\
        & 4.86E-02 $\pm$ 1.26E-02 & \texttt{repo} (2000)   & 6.15GB & 84 min & 9.7E-03 sec & 3.6E-03  \\
        \midrule
        \midrule
        \SetCell[r=3]{c} {KOL \\ with \texttt{iq4} kernel} 
        & Test Error & Dataset (size) & {CPU \\ Memory} & {Training \\ Time} & {Testing \\ Time} & {Pearson \\ (test)} \\
        \midrule
        &1.22E-02 & \SetCell[r=1]{c} \texttt{acti} (2000) & 1.33GB & 5 min & 5.6E-02 sec & 7.63E-05 \\
        &5.29E-03 & \texttt{repo} (2000)                  & 1.34GB & 5 min & 6.0E-02 sec & 7.64E-05 \\
        \bottomrule
    \end{tblr}
    \caption{Performance comparison of KOL on 3D unstructured datasets. Time single prediction test performed on a machine equipped with chip Apple M1 Pro.}
    \label{table:3D_unst_results}
\end{table}

\section{Conclusions}
\label{sec:conclusions}
In this study, we developed and formalized operator learning approaches to reconstruct activation and repolarization times in cardiac tissue, given an input activation region corresponding to electrically stimulated cells.
This problem is particularly significant for clinicians, as utilizing computational architectures for patient-specific simulations can improve clinical decision-making.
Although activation and repolarization times can be derived from PDE-based models, simulating these processes entails solving large-scale systems, resulting in high computational costs.
For this purpose, we adapted and evaluated two operator learning strategies: Fourier Neural Operators (FNO), based on the convolution theorem for the Fourier transform, and Kernel Operator Learning (KOL), based on kernel regression.
These trained operator learning techniques yield accurate and computationally efficient approximations of the target maps when evaluated on new samples.
Notably, in the repolarization case, we successfully approximated an operator map for which no corresponding PDE model is currently available.
Training data were generated by solving the monodomain model with spatial randomly distributed pulses with different ionic models using finite element method on 2D and 3D structured meshes, as well as on a physiologically realistic left ventricle.
Both methods demonstrated robust and accurate performance (with errors generally below 1\%) while significantly reducing computational costs compared to classical FEM-based simulations in a high number of evaluations.
Additionally, a systematic sensitivity analysis was conducted for the 2D case to assess hyperparameter dependence of both architectures.

Our numerical experiments demonstrate that KOL outperforms FNO in terms of accuracy and training time, even in the challenging case of 3D unstructured meshes. However, the computational gains of KOL are partially offset by the significant cost associated with kernel selection, which may pose a limitation.
This latter problem may be mitigated by choosing the kernel adaptively, \textit{e.g.} relying on the so called parametric Kernel Flows approaches \cite{owh19}.
Additionally, our results validate the feasibility of applying FNO to unstructured cardiac meshes, provided that suitable architectural modifications are implemented to accommodate non-uniform data structures.
Finally, while KOL offers superior accuracy, it comes at the expense of an increased inference time compared to the FNO counterpart. Hence, FNO can be the eligible choice for very high number of validation data.

\section*{Acknowledgments}
EC, NP, LP, SS, MV and GZ are members of INdAM-GNCS. EC and LFP have been supported by MUR (PRIN 202232A8AN\_002 and PRIN P2022B38NR\_001) funded by European Union - Next Generation EU. SS and GZ have been supported by MUR (PRIN 202232A8AN\_003 and PRIN P2022B38NR\_002) funded by European Union - Next Generation EU.

\newpage
\appendix
\setcounter{table}{0}
\renewcommand{\thetable}{\Alph{section}\arabic{table}}

\section{Extended Tables}
In this section we report the extended version of the tables presented in Sec.~\ref{sec:numerical}.

\begin{table}[H]
    \centering
    \small
    \begin{tblr}{
        colspec = {ccccccc}
    }
        \toprule
        \SetCell[r=12]{c} \rotatebox[origin=c]{90}{FNO} 
        & Test Error & Dataset (size) & lr policy & {GPU \\ Memory} & {Training \\ Time} \\
        \midrule
        & 3.33E-03 $\pm$ 8.41E-05 & \texttt{acti} (2000)     & \SetCell[r=6]{c} reduceOnPlateau & 5.62GB & 48 min \\
        & 2.66E-03 $\pm$ 2.27E-04 & \texttt{acti} (3000)     & & 5.79GB & 72 min \\
        & 3.33E-03 $\pm$ 3.13E-04 & \texttt{acti rot} (2000) & & 5.62GB & 48 min \\
        & 3.64E-03 $\pm$ 9.13E-04 & \texttt{repo} (2000)     & & 5.62GB & 48 min \\
        & 3.13E-03 $\pm$ 2.33E-04 & \texttt{repo} (3000)     & & 5.79GB & 72 min \\
        & 3.53E-03 $\pm$ 3.63E-04 & \texttt{repo rot} (2000) & & 5.62GB & 48 min \\
        \cmidrule{4}
        & 3.80E-03 $\pm$ 1.45E-04 & \texttt{acti} (2000)     & \SetCell[r=6]{c} None & 5.69GB & 48 min \\
        & 2.82E-03 $\pm$ 2.28E-04 & \texttt{acti} (3000)     & & 5.79GB & 72 min \\
        & 3.77E-03 $\pm$ 9.65E-04 & \texttt{acti rot} (2000) & & 5.62GB & 48 min \\
        & 3.40E-03 $\pm$ 5.77E-05 & \texttt{repo} (2000)     & & 5.62GB & 48 min \\
        & 3.53E-03 $\pm$ 3.32E-04 & \texttt{repo} (3000)     & & 5.79GB & 72 min \\
        & 4.16E-03 $\pm$ 4.04E-05 & \texttt{repo rot} (2000) & & 5.62GB & 48 min \\
        \bottomrule
    \end{tblr}
    \caption{Performance comparison of FNO on 2D datasets.}
    \label{table:2D_results_FNO_appendix}
\end{table}

\begin{table}[H]
    \centering
    \small
    \begin{tblr}{
        colspec = {ccccccc}
    }
        \toprule
        \SetCell[r=34]{c} \rotatebox[origin=c]{90}{KOL} 
        & Test Error & Dataset (size) & Kernel (iq) & {CPU \\ Memory} & {Training\\ Time} \\
        \midrule
        & 9.33E-04 & \texttt{acti rot} (2000) & iq-1& 0.91GB & 476 sec\\
        & 9.35E-04 & \texttt{acti rot} (2000) & iq-2& 0.92GB & 517 sec\\
        & 9.39E-04 & \texttt{acti rot} (2000) & iq-3& 0.91GB & 517 sec\\
        & 9.34E-04 & \texttt{acti rot} (2000) & iq-4& 0.91GB & 459 sec\\
        & 1.15E-03 & \texttt{acti rot} (2000) & iq-5& 0.97GB & 477 sec\\
        & 1.19E-01 & \texttt{acti rot} (2000) & ntk-1& 0.97GB & 660 sec\\
        & 1.41E-01 & \texttt{acti rot} (2000) & ntk-2& 0.97GB & 691 sec\\
        & 1.42E-01 & \texttt{acti rot} (2000) & ntk-3& 0.97GB & 687 sec\\
        & 1.73E-01 & \texttt{acti rot} (2000) & rbf-1& 0.92GB & 93 sec\\
        & 1.19E-01 & \texttt{acti rot} (2000) & rbf-2& 0.92GB & 92 sec\\
        & 1.18E-01 & \texttt{acti rot} (2000) & rbf-3& 0.92GB & 89 sec\\
        \cmidrule{3}
        & 1.07E-03 & \texttt{acti} (2000) & iq-1 & 0.91GB & 444 sec\\
        & 1.08E-03 & \texttt{acti} (2000) & iq-2 & 0.91GB & 451 sec\\
        & 1.09E-03 & \texttt{acti} (2000) & iq-3 & 0.92GB & 442 sec\\
        & 1.08E-03 & \texttt{acti} (2000) & iq-4 & 0.92GB & 448 sec\\
        & 1.38E-03 & \texttt{acti} (2000) & iq-5 & 0.91GB & 448 sec\\
        & 1.31E-01 & \texttt{acti} (2000) & ntk-1 & 0.97GB & 656 sec\\
        & 1.54E-01 & \texttt{acti} (2000) & ntk-2 & 0.97GB & 652 sec\\
        & 1.55E-01 & \texttt{acti} (2000) & ntk-3 & 0.97GB & 636 sec\\
        & 1.86E-01 & \texttt{acti} (2000) & rbf-1 & 0.92GB & 82 sec\\
        & 1.31E-01 & \texttt{acti} (2000) & rbf-2 & 0.92GB & 80 sec\\
        & 1.30E-01 & \texttt{acti} (2000) & rbf-3 & 0.92GB & 80 sec\\
        \cmidrule{3}
        & 9.51E-04 & \texttt{acti} (3000) & iq-1 & 1.15GB & 613 sec\\
        & 9.53E-04 & \texttt{acti} (3000) & iq-2 & 1.15GB & 609 sec\\
        & 9.54E-04 & \texttt{acti} (3000) & iq-3 & 1.15GB & 616 sec\\
        & 9.52E-04 & \texttt{acti} (3000) & iq-4 & 1.15GB & 611 sec\\
        & 1.11E-03 & \texttt{acti} (3000) & iq-5 & 1.15GB & 611 sec\\
        & 1.28E-01 & \texttt{acti} (3000) & ntk-1 & 1.20GB & 979 sec\\
        & 1.52E-01 & \texttt{acti} (3000) & ntk-2 & 1.20GB & 980 sec\\
        & 1.51E-01 & \texttt{acti} (3000) & ntk-3 & 1.20GB & 1007 sec\\
        & 1.83E-01 & \texttt{acti} (3000) & rbf-1 & 1.18GB & 117 sec\\
        & 1.28E-01 & \texttt{acti} (3000) & rbf-2 & 1.15GB & 125 sec\\
        & 1.27E-01 & \texttt{acti} (3000) & rbf-3 & 1.15GB & 127 sec\\

        \bottomrule
   \end{tblr}
    \caption{Performance comparison of KOL on 2D datasets for reconstructing activation maps.}
    \label{table:2D_results_kol_acti_appendix}
\end{table}

\begin{table}[H]
    \centering
    \small
    \begin{tblr}{
        colspec = {ccccccc}
    }
        \toprule
        \SetCell[r=34]{c} \rotatebox[origin=c]{90}{KOL} 
        & Test Error & Dataset (size) & Kernel (iq) & {CPU \\ Memory} & {Training\\ Time} \\
        \midrule
        & 4.69E-04 & \texttt{repo rot} (2000) & iq-1& 0.91GB & 418 sec\\
        & 4.67E-04 & \texttt{repo rot} (2000) & iq-2& 0.91GB & 389 sec\\
        & 4.69E-04 & \texttt{repo rot} (2000) & iq-3& 0.92GB & 389 sec\\
        & 4.69E-04 & \texttt{repo rot} (2000) & iq-4& 0.91GB & 393 sec\\
        & 6.72E-04 & \texttt{repo rot} (2000) & iq-5& 0.91GB & 397 sec\\
        & 5.95E-02 & \texttt{repo rot} (2000) & ntk-1& 0.97GB & 575 sec\\
        & 7.02E-02 & \texttt{repo rot} (2000) & ntk-2& 0.97GB & 571 sec\\
        & 6.98E-02 & \texttt{repo rot} (2000) & ntk-3& 0.97GB & 565 sec\\
        & 8.33E-02 & \texttt{repo rot} (2000) & rbf-1& 0.92GB & 72 sec\\
        & 5.93E-02 & \texttt{repo rot} (2000) & rbf-2& 0.92GB & 71 sec\\
        & 5.91E-02 & \texttt{repo rot} (2000) & rbf-3& 0.92GB & 72 sec\\
        \cmidrule{3}
        & 4.93E-04 & \texttt{repo} (2000) & iq-1 & 0.91GB & 408 sec\\
        & 4.92E-04 & \texttt{repo} (2000) & iq-2 & 0.91GB & 450 sec\\
        & 5.19E-04 & \texttt{repo} (2000) & iq-3 & 0.91GB & 473 sec\\
        & 4.93E-04 & \texttt{repo} (2000) & iq-4 & 0.91GB & 473 sec\\
        & 7.82E-04 & \texttt{repo} (2000) & iq-5 & 0.92GB & 414 sec\\
        & 6.49E-02 & \texttt{repo} (2000) & ntk-1 & 0.97GB & 595 sec\\
        & 7.64E-02 & \texttt{repo} (2000) & ntk-2 & 0.97GB & 632 sec\\
        & 7.61E-02 & \texttt{repo} (2000) & ntk-3 & 0.97GB & 622 sec\\
        & 9.01E-02 & \texttt{repo} (2000) & rbf-1 & 0.92GB & 74 sec\\
        & 6.47E-02 & \texttt{repo} (2000) & rbf-2 & 0.92GB & 76 sec\\
        & 6.44E-02 & \texttt{repo} (2000) & rbf-3 & 0.92GB & 90 sec\\
        \cmidrule{3}
        & 4.74E-04 & \texttt{repo} (3000) & iq-1 & 1.15GB & 592 sec\\
        & 4.72E-04 & \texttt{repo} (3000) & iq-2 & 1.15GB & 576 sec\\
        & 4.97E-04 & \texttt{repo} (3000) & iq-3 & 1.15GB & 574 sec\\
        & 4.74E-04 & \texttt{repo} (3000) & iq-4 & 1.15GB & 582 sec\\
        & 7.19E-04 & \texttt{repo} (3000) & iq-5 & 1.15GB & 573 sec\\
        & 6.20E-02 & \texttt{repo} (3000) & ntk-1 & 1.20GB & 1064 sec\\
        & 7.39E-02 & \texttt{repo} (3000) & ntk-2 & 1.20GB & 1007 sec\\
        & 7.34E-02 & \texttt{repo} (3000) & ntk-3 & 1.20GB & 987 sec\\
        & 8.77E-02 & \texttt{repo} (3000) & rbf-1 & 1.15GB & 120 sec\\
        & 6.19E-02 & \texttt{repo} (3000) & rbf-2 & 1.15GB & 117 sec\\
        & 6.16E-02 & \texttt{repo} (3000) & rbf-3 & 1.15GB & 117 sec\\

        \bottomrule
   \end{tblr}
    \caption{Performance comparison of KOL on 2D datasets for reconstructing repolarization maps.}
    \label{table:2D_results_kol_repo_appendix}
\end{table}

\begin{table}[H]
    \centering
    \small
    \begin{adjustbox}{max width=\textwidth}
    \begin{tblr}{
        colspec = {ccccccc}
    }
        \toprule
        Test Error & Dataset (size) & L & width & parameters & {GPU \\ Memory} & {Training \\ Time}  & {Testing \\ Time} & {Pearson \\ (test)} \\
        \midrule
        1.49E-01 $\pm$ 3.43E-02 & \SetCell[r=2]{c} \texttt{acti} (1000)  & 3 & 32 & 10.9M & 6.06GB & 58 min  & 9.4E-03 sec & 2.3E-02 \\
        1.27E-01 $\pm$ 1.46E-03 &                               & 4 & 32 & 10.9M & 6.07GB & 59 min  & 1.0E-02 sec & 2.6E-02 \\
        \cmidrule{2}
        1.44E-01 $\pm$ 2.30E-02 & \SetCell[r=9]{c} \texttt{acti} (2000)  & 1 & 2   & 10.8M & 6.12GB & 81 min & 8.6E-03 sec & 2.6E-02 \\
        7.93E-02 $\pm$ 6.44E-03 &                               & 1 & 4   & 10.8M & 6.12GB & 80 min & 8.3E-03 sec & 9.7E-03 \\
        7.63E-02 $\pm$ 1.97E-03 &                               & 1 & 8   & 10.8M & 6.12GB & 81 min & 8.3E-03 sec & 9.2E-03 \\
        7.18E-02 $\pm$ 4.03E-03 &                               & 1 & 16  & 10.9M & 6.13GB & 79 min & 8.4E-03 sec & 6.8E-03 \\
        8.81E-02 $\pm$ 2.56E-02 &                               & 3 & 32  & 11.0M & 6.15GB & 83 min & 9.5E-03 sec & 2.3E-02 \\
        1.14E-01 $\pm$ 1.72E-02 &                               & 2 & 32  & 10.9M & 6.14GB & 82 min & 8.9E-03 sec & 1.1E-02 \\
        7.15E-02 $\pm$ 9.84E-03 &                               & 1 & 32  & 10.9M & 6.13GB & 80 min & 8.5E-03 sec & 4.3E-03 \\
        1.07E-01 $\pm$ 2.57E-02 &                               & 1 & 64  & 11.0M & 6.14GB & 80 min & 8.5E-03 sec & 2.8E-02 \\
        1.49E-01 $\pm$ 1.12E-02 &                               & 1 & 128 & 11.3M & 6.17GB & 81 min & 9.4E-03 sec & 3.3E-02 \\
        \cmidrule{2}
        4.86E-02 $\pm$ 1.26E-02 & \texttt{repo} (2000)                   & 3 & 32  & 10.9M & 6.15GB & 84 min & 9.7E-03 sec & 3.6E-03  \\
        \bottomrule
    \end{tblr}
    \end{adjustbox}
    \caption{Performance comparison of FNO on 3D unstructured datasets. Time single prediction test performed on a machine equipped with chip Apple M1 Pro.}
    \label{table:3D_unst_results_appendix}
\end{table}

\section{Nomenclature of kernels for KOL}
\label{app:nomenKOL}
\begin{table}[H]
    \centering
    \small
    \begin{tblr}{
        colspec = {ccccccc}
    }
        \toprule
        {Name} & {Kernel\\Type}& {$\sigma$} &  $d_{\mathrm{nn}}$ & {Activation\\Function} & {$\sigma_1$} & {$\sigma_2$}  \\
        \midrule
        \texttt{iq1} & IQ & / & / & / & 1E-5 & 1E-2\\
        \cmidrule{1}
        \texttt{iq2} & IQ & / & / & / & 1E-5 & 1E-1 \\    
        \cmidrule{1}
        \texttt{iq3} & IQ & / & / & / & 1E-4 & 1E-2\\
        \cmidrule{1}
        \texttt{iq4} & IQ & / & / & / & 1E-4 & 1E-1\\
        \cmidrule{1}
        \texttt{iq5} & IQ & / & / & / & 1E-3 & 1E-2\\
        \cmidrule{1}
        \texttt{rbf1} & RBF & 1 & / & / & / & /\\
        \cmidrule{1}
        \texttt{rbf2} & RBF & 10 & / & / & / & / \\
        \cmidrule{1}
        \texttt{rbf3} & RBF & 100 & / & / & / & /\\
        \cmidrule{1}
        \texttt{ntk1} & NTK & / & 3 & Sigmoid & / & / \\
        \cmidrule{1}
        \texttt{ntk2} & NTK & / & 4 & Sigmoid & / & / \\
        \cmidrule{1}
        \texttt{ntk3} & NTK & / & 3 & ReLu & / & / \\
        \bottomrule
    \end{tblr}
    \caption{Nomenclature of kernel functions tested for KOL.}
    \label{table:nomenKOL_appendix}
\end{table}
\section{Repolarization error distribution: 3D slab}
\begin{figure}[H]
\centering
\begin{subfigure}{0.24\textwidth}
    \centering
    \includegraphics[width=\textwidth]{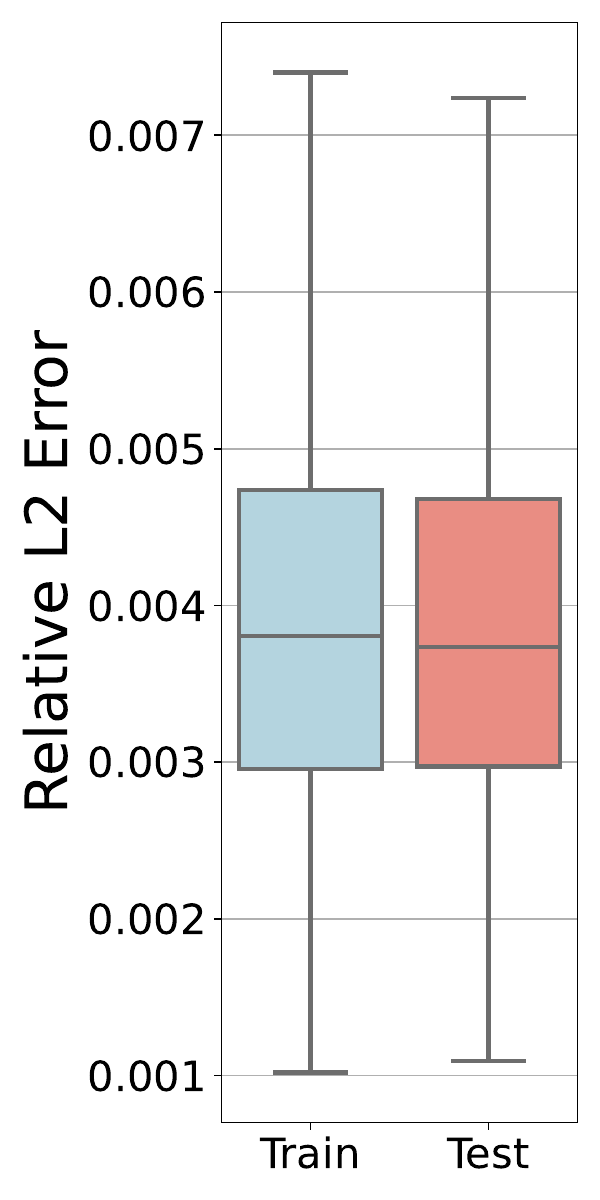}
    \caption{FNO box plot.}
    \label{fig:fno_box_repo}
\end{subfigure}
    \begin{subfigure}{0.55\textwidth}
    \centering
    \includegraphics[width=\textwidth]{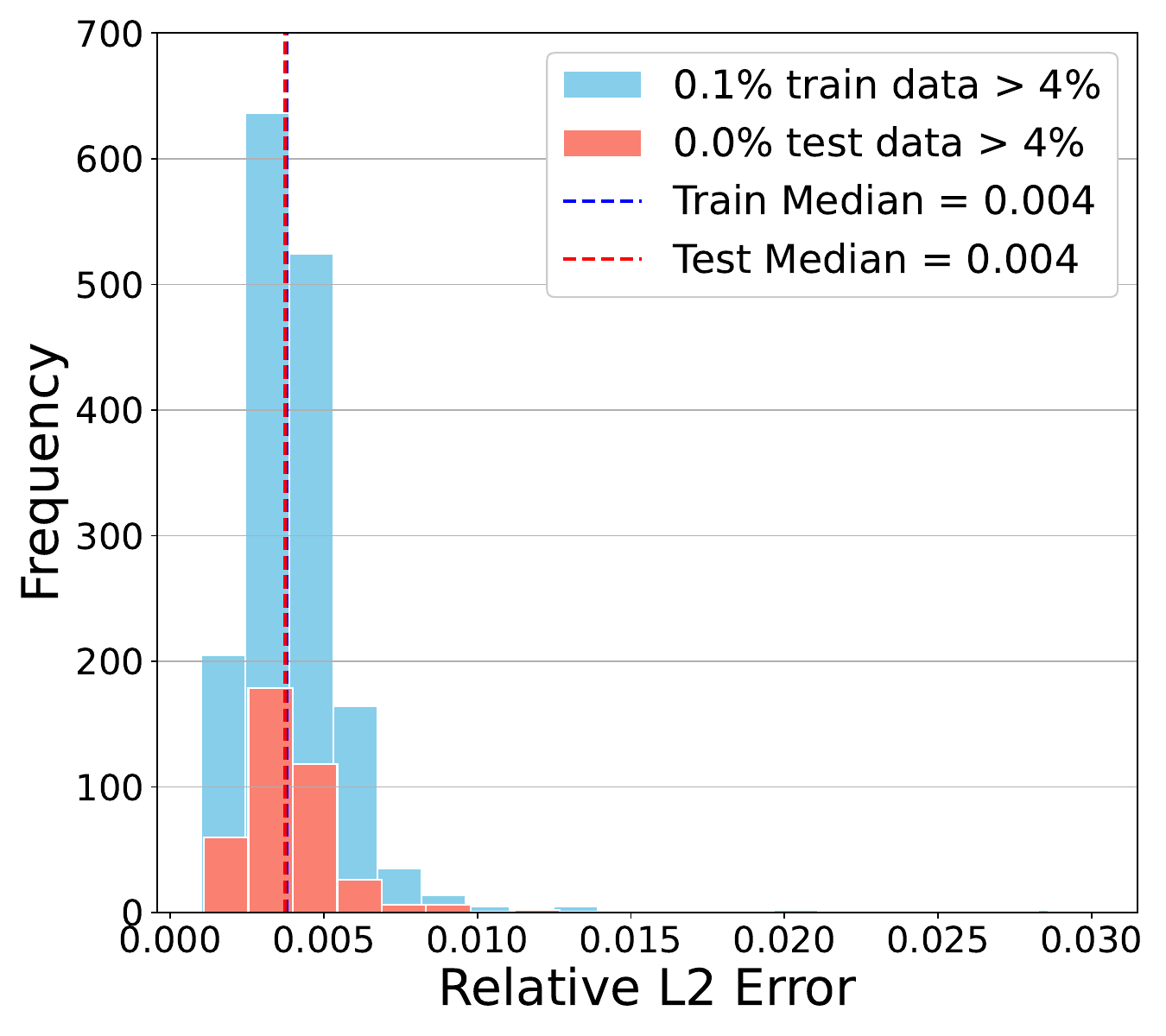}
    \caption{FNO histogram.}
    \label{fig:fno_hist_repo}
\end{subfigure}
\caption{FNO box plot and histogram for 3D dataset \texttt{repo} 2000 relative to the best model trained (outliers not shown). For the training set, 0.1\% of the data have a relative L2 error greater than 4\%, while for the test set no data exceed this threshold.}
\label{fig:statsFNOrepo}
\end{figure}

\begin{figure}[H]
\centering
\begin{subfigure}{0.18\textwidth}
    \centering
    \includegraphics[width=\textwidth]{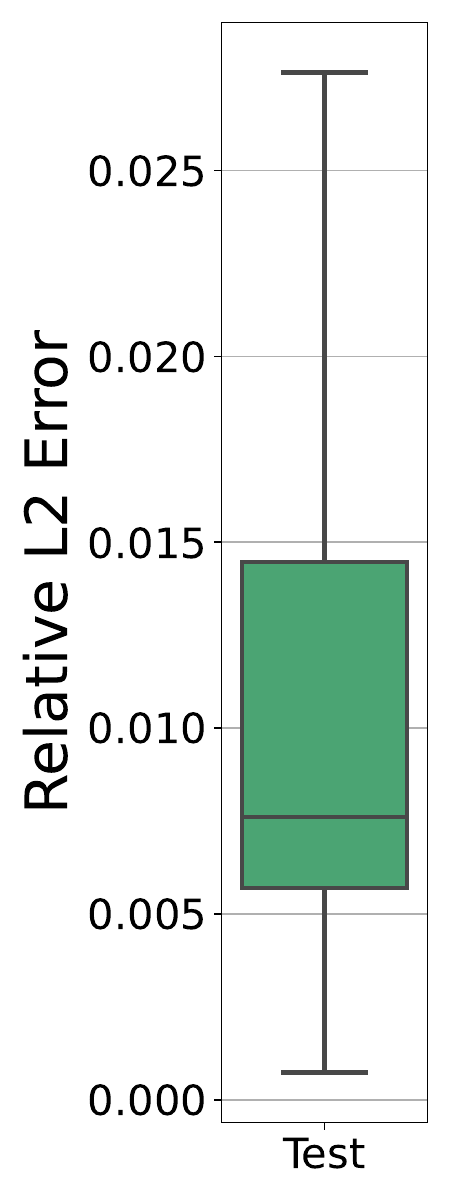}
    \caption{KOL box plot.}
    \label{fig:kol_box_repo}
\end{subfigure}
    \begin{subfigure}{0.55\textwidth}
    \centering
    \includegraphics[width=\textwidth]{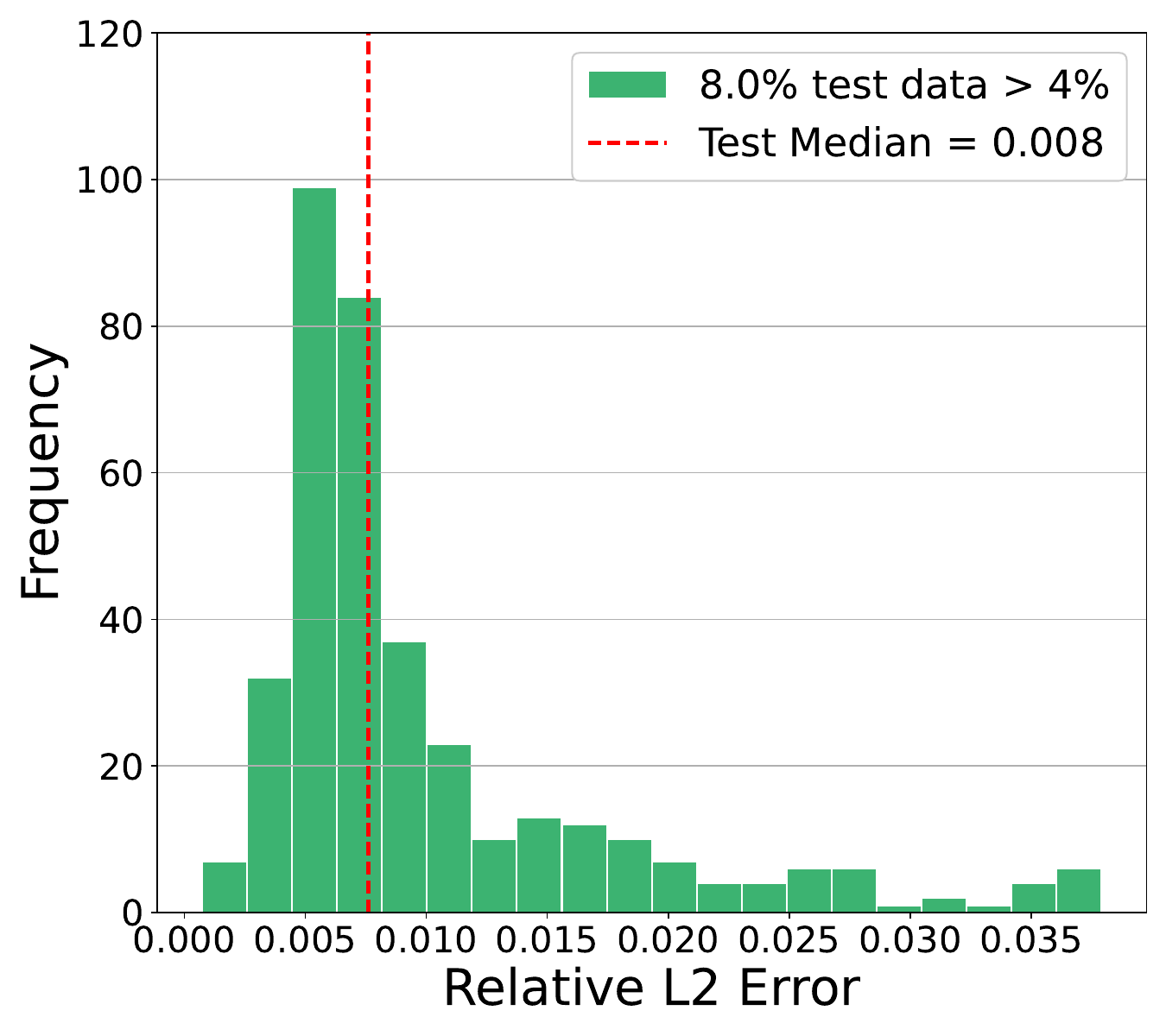}
    \caption{KOL histogram.}
    \label{fig:kol_hist_repo}
\end{subfigure}
\caption{KOL box plot and histogram for 3D dataset \texttt{repo} 2000 relative to the best model trained (outliers not shown). Training results are not shown since we achieve machine precision. For the test set, 8\% of the data have a relative L2 error greater than 4\%.}
\label{fig:statsKOLrepo}
\end{figure}

\end{document}